\documentclass[12pt]{article}
\usepackage[dvipdfmx]{graphicx}
\usepackage{color}
\usepackage[all]{xy}

\usepackage{here}
\usepackage{amsmath,amssymb,amsthm, amscd}
\usepackage{ytableau,mathrsfs}
  
  \makeatletter
  \@addtoreset{equation}{section}
  \makeatother

\renewcommand{\tilde}[1]{\widetilde{#1}}
\newcommand\floor[1]{\lfloor#1\rfloor}

\newcommand{\jump}[1]{\ensuremath{[\![#1]\!]} }
\newtheorem{dfn}{Definition}[section]
\newtheorem{thm}[dfn]{Theorem}
\newtheorem{lem}[dfn]{Lemma}
\newtheorem{prop}[dfn]{Proposition}

\theoremstyle{remark}
\newtheorem{rem}{Remark}[section]

\begin{document}
\begin{center} 
{\bf {\LARGE Symplectic resolutions of the Hilbert squares of ADE surface singularities}}
\end{center} 
\vspace{0.4cm}

\begin{center}
{\large Ryo Yamagishi}
\end{center} 
\vspace{0.4cm}

\begin{abstract}
We study symplectic resolutions of the Hilbert scheme of two points on a surface with one ADE-singularity. We also characterize such singularities by central fibers of their symplectic resolutions. As an application, we show that these singularities are isomorphic to the Slodowy slices which are transversal to the `sub-subregular' orbits in the nilpotent cones of ADE-types.
\end{abstract}

\section{Introduction}\label{1}
Throughout this paper we work over the field of complex numbers.

The concept of symplectic singularities, first introduced in \cite{Bea}, is a singular version of symplectic manifolds and have been studied by many authors. The only examples of two-dimensional symplectic singularities are ADE-singularities. However,  four-dimensional (and higher) ones are not classified. One strategy for the classification problem may be to focus on the central fiber of a symplecic resolution of a given symplectic singularity. In the two-dimensional cases, the singularities are indeed determined by the central fibers, which are trees of projective lines. In this paper we give more evidence for this strategy by showing that certain four-dimensional symplectic singularities are determined by the central fibers of their symlectic resolutions.

Each of the singularities that we will treat is obtained as the Hilbert scheme $\mathrm{Hilb}^2(\Gamma)$ of two points on a surface $\Gamma$ with one ADE-singularity. These singularities are indeed symplectic singularities and will turn out to have unique symplectic resolutions (see Section \ref{2}). The resolutions can be obtained by just successively blowing up the singular loci, but it is not easy to calculate the central fibers directly. Instead, we consider a symplectic resolution $\mathrm{Hilb}^2(\tilde{\Gamma})$ of the symmetric product $\mathrm{Sym}^2(\Gamma)$ and perform Mukai flops to obtain a symplectic resolution $\mathcal{H}$ of $\mathrm{Hilb}^2(\Gamma)$ where $\tilde{\Gamma}$ is the minimal resolution of $\Gamma$.

After determining the central fiber $F$ of $\mathcal{H}$, we study the formal neighborhood of the fiber. In the two-dimensional cases, it is well-known that the formal neighborhood of a smooth rational curve $C$ in a smooth symplectic surface is uniquely determined as (the completion of) the cotangent bundle $T_C^*$. Also in the case of $\mathrm{Hilb}^2(\Gamma)$, the symplectic structure and the induced Poisson structure gives strong restriction on the formal neighborhood of $F$. We will show that the formal neighborhood of $F$ is determined by $F$ itself and give explicit description of it in Section \ref{3}. As a consequence, we can show that the analytic type of the singularity is also determined by the isomorphism class of the central fiber (Theorem \ref{thm:main}).

We also consider certain symplectic varieties called {\em Slodowy slices} in complex simple Lie algebras $\mathfrak{g}$ of ADE-types. Slodowy slices give many examples of symplectic varieties and have been studied also from the representation-theoretic viewpoint. The well-known result of Brieskorn-Slodowy theory states that the Slodowy slices which are transversal to the subregular orbits in the nilpotent cones have ADE-singularities of the same types as $\mathfrak{g}$.

In this paper we particularly treat the Slodowy slices which are transversal to the `sub-subregular' orbits (which mean the third biggest orbits) in the nilpotent cones of ADE-types. For types $\mathbf{A}_n\,(n\ge3)$ and $\mathbf{E}_n$, the sub-subregular orbits are unique and the slices have unique `bad' points which are called 0-dimensional symplectic leaves (see Theorem \ref{thm:Kaledin}). Similarly, $\mathrm{Hilb}^2(\Gamma)$ have unique 0-dimensional symplectic leaves for types $\mathbf{A}_n\,(n\ge3)$ and $\mathbf{E}_n$. For type $\mathbf{D}_n$, the sub-subregular orbits are not unique. (There are three if $n=4$ and two if $n\ge5$.) However, the 0-dimensional symplectic leaves of $\mathrm{Hilb}^2(\Gamma)$ are also non-unique, and there is still a bijection between the 0-dimensional leaves of $\mathrm{Hilb}^2(\Gamma)$ and the slice for $\mathbf{D}_n$.

The fibers of symplectic resolutions of Slodowy slices are known as Springer fibers. In particular, the sub-subregular cases are studied by Lorist \cite{Lo}. These fibers are in fact isomorphic to the ones for $\mathrm{Hilb}^2(\Gamma)$. Thus, our result above can apply to show that the Slodowy slices have the same singularity types as $\mathrm{Hilb}^2(\Gamma)$ (Theorem \ref{thm:slice}). Therefore, this result can be regarded as an analogue of the classical result for subregular orbits mentioned above.

This theorem (Theorem \ref{thm:slice}) for type $\mathbf{A}_n$ also follows from the work by Manolescu \cite[Thm. 1.1]{M}. More strongly, his result shows that, for each $k,n\in\mathbb{N}$ with $2k\le n+1$, there is a Slodowy slice in the nilpotent cone of $\mathfrak{sl}_{n+1}$ that has an open immersion to the Hilbert scheme $\mathrm{Hilb}^k(\Gamma)$ with $\Gamma$ of type $\mathbf{A}_n$. Then, our natural next task should be to extend his result to $\mathbf{D}_n$ and $\mathbf{E}_n$. Moreover, it also seems interesting to find analogous results for non-simply laced types of Lie algebras like the Brieskorn-Slodowy theory for subregular orbits.

The contents of this article are summarized as follows. In Section \ref{2}, we study $\mathrm{Hilb}^2(\Gamma)$ and its symplectic resolution $\mathcal{H}$. In particular we give explicit defining equations of $\mathrm{Hilb}^2(\Gamma)$ and explicit procedures of Mukai flops to obtain $\mathcal{H}$ from $\mathrm{Hilb}^2(\tilde{\Gamma})$. Also, we describe the central fibers $F$ in the resolution $\mathcal{H}$. In Section \ref{3}, we study the formal neighborhoods of the fibers $F$ using the symplectic and Poisson structures. We will show that the formal neighborhood of a Lagrangian submanifold is locally a cotangent bundle (Proposition \ref{prop:local}) and that the global structure is controlled by some cohomology groups (Lemma \ref{lem:normal}). Then we will see that the isomorphism type of the neighborhood of each irreducible component of $F$ is uniquely determined. Next we will show that the way of gluing of such neighborhoods is also unique up to isomorphism, and thereby prove the main theorem (Theorem \ref{thm:main}). In Section \ref{4}, we introduce the Slodowy slices and apply the results in Section \ref{3} to show the coincidence of the singularities of the Slodowy slices with those of $\mathrm{Hilb}^2(\Gamma)$ (Theorem \ref{thm:slice}).

\vspace{5mm}
\noindent{\bf Acknowledgements}\\
The author would like to thank Yoshinori Namikawa for many helpful comments. The author is supported by JSPS KAKENHI Grant Number JP16J04485.

\section{$\mathrm{Hilb}^2(\Gamma)$ for a surface $\Gamma$ with an ADE-singularity}\label{2}
In this section we investigate properties of the Hilbert scheme of two points on a surface $\Gamma$ with one singular point of type $T=\mathbf{A}_n,\mathbf{D}_n$ or $\mathbf{E}_n$. We may assume that $\Gamma$ is defined as the specific hypersurface $\Gamma=\mathrm{Spec}\,\mathbb{C}[x,y,z]/(f)\subset\mathbb{C}^3$ where
{\footnotesize
$$f=\begin{cases}x^{n+1}-yz&\hspace{5mm}\text{if }\;T=\mathbf{A}_n\\
x^{n-1}+xy^2+z^2&\hspace{5mm}\text{if }\;T=\mathbf{D}_n\\
x^4+y^3+z^2&\hspace{5mm}\text{if }\;T=\mathbf{E}_6\\
x^3y+y^3+z^2&\hspace{5mm}\text{if }\;T=\mathbf{E}_7\\
x^5+y^3+z^2&\hspace{5mm}\text{if }\;T=\mathbf{E}_8\\
\end{cases}$$}
since we are only interested in a (formal) neighborhood of the subvariety of the Hilbert scheme consisting of points supported at the singular point. Let $\mathrm{Sym}^2(\Gamma)$ be the symmetric product of two copies of $\Gamma$. Then we have a birational morphism $\phi:\mathrm{Hilb}^2(\Gamma)\to\mathrm{Sym}^2(\Gamma)$ which is called the Hibert-Chow morphism. $\phi$ takes a 0-dimensional subscheme of $\Gamma$ to its support. See \cite[\S2]{Fo} or \cite[Thm. 2.7]{Ber} for the precise construction.

We also have the Hilbert-Chow morphism $\phi':\mathrm{Hilb}^2(\mathbb{C}^3)\to\mathrm{Sym}^2(\mathbb{C}^3)$ for $\mathbb{C}^3$. It is known that $\phi'$ is the blowing-up of the diagonal and that $\phi'|_{\mathrm{Hilb}^2(\Gamma)}=\phi$ (see e.g. \cite[Ex. 2.7 and Rem. 2.9]{Ber}). Since $\mathrm{Hilb}^2(\Gamma)$ is irreducible \cite[Thm. A]{Z}, we have the following result using the property of a blowing-up.

\begin{prop}\label{prop:HC}
The Hilbert-Chow morphism $\phi:\mathrm{Hilb}^2(\Gamma)\to\mathrm{Sym}^2(\Gamma)$ can be identified with the blowing-up of $\mathrm{Sym}^2(\Gamma)$ along the reduced diagonal.
\end{prop}

\vspace{3mm}

The same result is in fact true for $\mathrm{Hilb}^d(\Gamma)$ with any $d\in\mathbb{N}$ when $T=\mathbf{A}_n$ \cite[Thm. D]{Z}. In the subsequent subsections, we will calculate the defining equation of $\mathrm{Hilb}^2(\Gamma)$ for each $T$ by using Proposition \ref{prop:HC}. 

\subsection{Case $T=\mathbf{A}_n$}\label{2.1}
In this subsection we consider the case when $T=\mathbf{A}_n$. Then the surface $\Gamma_n:=\Gamma$ is the quotient of $\mathbb{C}^2$ by the action of the cyclic group $C_{n+1}$ generated by {\footnotesize$\begin{pmatrix}\zeta&0\\0&\zeta^{-1}\end{pmatrix}$} where $\zeta\in\mathbb{C}$ is a primitive $(n+1)$-th root of unity. Therefore, $\mathrm{Sym}^2(\Gamma_n)$ is the quotient of $\mathrm{Spec}\,(\mathbb{C}[a_1,a_2]\otimes\mathbb{C}[b_1,b_2])=\mathbb{C}^2\times\mathbb{C}^2$ by the natural action of the wreath product $G=C_{n+1}\wr \mathfrak{S}_2$.
The ring of invariants by the normal subgroup $(C_{n+1})^{\times 2}\lhd G$ is generated by
$$u_1=a_1^{n+1},\;v_1=a_1a_2,\;w_1=a_2^{n+1},\;u_2=b_1^{n+1},\;v_2=b_1b_2,\;w_2=b_2^{n+1}.$$
Since the residual action of $\mathfrak{S}_2$ swaps $a_i$ and $b_i$, the ring of invariants by $G$ is generated by the following elements:
{\footnotesize
$$\begin{aligned}x_1=u_1+u_2,\;
x_2&=v_1+v_2,\;
x_3=w_1+w_2,\\
x_4=(u_1-u_2)^2,\;
x_5&=(v_1-v_2)^2,\;
x_6=(w_1-w_2)^2,\\
x_7=(u_1-u_2)(v_1-v_2),\;
x_8&=(u_1-u_2)(w_1-w_2),\;
x_9=(v_1-v_2)(w_1-w_2).
\end{aligned}$$}

\begin{prop}\label{prop:relation1}
The ideal of relations of $x_i$'s is generated by the following ten elements of $\mathbb{C}[x_1,x_2,\dots,x_9]$:
{\footnotesize
$$\begin{aligned}f_1&=\sum_{i=0}^{\floor{(n+1)/2}}
\begin{pmatrix}n+1\\2i\end{pmatrix}x_2^{n-2i+1}x_5^i-2^{n-1}(x_1x_3+x_8),\\
f_2&=\sum_{i=0}^{\floor{(n+1)/2}}
\begin{pmatrix}n+1\\2i+1\end{pmatrix}x_2^{n-2i}x_5^{i+1}-2^{n-1}(x_1x_9+x_3x_7),\\
f_3&=x_5x_8-x_7x_9,\;f_4=x_5x_6-x_9^2,\;f_5=x_4x_5-x_7^2,\\
f_6&=\sum_{i=0}^{\floor{(n+1)/2}}
\begin{pmatrix}n+1\\2i+1\end{pmatrix}x_2^{n-2i}x_5^i x_9-2^{n-1}(x_1x_6+x_3x_8),\\
f_7&=\sum_{i=0}^{\floor{(n+1)/2}}
\begin{pmatrix}n+1\\2i+1\end{pmatrix}x_2^{n-2i}x_5^i x_7-2^{n-1}(x_3x_4+x_1x_8),\\
f_8&=x_7x_8-x_4x_9,\,f_9=x_6x_7-x_8x_9,\;f_{10}=x_4x_6-x_8^2.
\end{aligned}$$}
\end{prop}

{\em Proof.} We should show that the kernel of the composite of the ring maps $f:\mathbb{C}[x_1,x_2,\dots,x_9]\to\mathbb{C}[u_1,v_1,w_1,u_2,v_2,w_2]$ and $g:\mathbb{C}[u_1,v_1,w_1,u_2,v_2,w_2]\to\mathbb{C}[a_1,a_2,b_1,b_2]$ is generated by $f_i$'s. Note that $\mathrm{Ker}\,g$ is generated by $r_i=u_i w_i-v_i^{n+1},\,i=1,2$. The $\mathfrak{S}_2$-action on $\mathbb{C}[a_1,a_2,b_1,b_2]$ naturally induces an action on $\mathbb{C}[u_1,v_1,w_1,u_2,v_2,w_2]$ which switches $r_1$ and $r_2$, and its invariant subring is equal to $\mathrm{Im}\,f$. Therefore, we see that $\mathrm{Ker}\,g\cap\mathrm{Im}\,f$, as an ideal of $\mathrm{Im}\,f$, is generated by $r_1+r_2,(v_1-v_2)(r_1-r_2),(w_1-w_2)(r_1-r_2)$ and $(u_1-u_2)(r_1-r_2)$. One can check that $f_1,f_2,f_6$ and $f_7$ respectively are mapped by $f$ to these generators up to constants using the following equalities:
{\footnotesize
$$\sum_{i=0}^{\floor{(n+1)/2}}
\begin{pmatrix}n+1\\2i\end{pmatrix}x_2^{n-2i+1}x_5^i=\frac{1}{2}\{(x_2+\sqrt{x_5})^{n+1}+(x_2-\sqrt{x_5})^{n+1}\},$$
$$\sum_{i=0}^{\floor{(n+1)/2}}
\begin{pmatrix}n+1\\2i+1\end{pmatrix}x_2^{n-2i}x_5^i=\frac{1}{2\sqrt{x_5}}\{(x_2+\sqrt{x_5})^{n+1}-(x_2-\sqrt{x_5})^{n+1}\}.$$}
The claim follows since $\mathrm{Ker}\,f$ is clearly generated by the other six $f_i$'s.
\qed

\vspace{3mm}

Now we can explicitly calculate the blowing-up $\mathrm{Bl}_\Delta$ of $\mathrm{Sym}^2(\Gamma_n)$ along the diagonal $\Delta=\{x_4=x_5=\cdots=x_9=0\}$. Let $y_4,\dots,y_9$ be the homogeneous coordinates of $\mathrm{Bl}_\Delta$ corresponding to the coordinates $x_4,\dots,x_9$. Then the fiber of the origin is a Veronese surface in $\mathbb{P}^5=\mathrm{Proj}\,\mathbb{C}[y_4,\dots,y_9]$ and hence isomorphic to $\mathbb{P}^2$. Let $U_i=\{y_i\ne0\}\subset \mathrm{Bl}_\Delta\,(i=4,\dots,9)$ be the affine open covering. Then except the case $i=5$ and $n\ge3$, the $U_i$'s have hypersurface singularities whose singular loci are smooth. For $i=5$ and $n\ge3$, we set $x_i'=y_i/y_5,\;i=4,\dots,9$. Then we can eliminate variables, and $U_5$ is defined by the following two equations in $\mathbb{C}^6=\mathrm{Spec}\,\mathbb{C}[x_1,x_2,x_3,x_5,x_7',x_9']$:
{\footnotesize
$$\begin{aligned}\tilde{f}_1&=\sum_{i=0}^{\floor{(n+1)/2}}
\begin{pmatrix}n+1\\2i\end{pmatrix}x_2^{n-2i+1}x_5^i-2^{n-1}(x_1x_3+x_5x_7'x_9')\\
\tilde{f}_2&=\sum_{i=0}^{\floor{(n+1)/2}}
\begin{pmatrix}n+1\\2i+1\end{pmatrix}x_2^{n-2i}x_5^i-2^{n-1}(x_1x_9'+x_3x_7').
\end{aligned}$$}

The (reduced) singular locus of $U_5$ is defined as
$$\{x_1-x_2x_7'=x_3-x_2x_9'=x_5-x_2^2=x_2^{n-1}-x_7'x_9'=0\}.$$ Thus, it has an $\mathbf{A}_{n-2}$-singular point at the origin.

To describe properties of $\mathrm{Hilb}^2(\Gamma_n)$, we introduce symplectic varieties.

\begin{dfn}
$\boldsymbol{\cdot}$ (cf. \cite{Bea}) A normal (singular) variety $X$ is a {\bf symplectic variety} (or has {\bf symplectic singularities}) if the smooth part $X_\mathrm{reg}$ of $X$ has a holomorphic symplectic form $\omega$ and if there is a desingularization $\pi:Y\to X$ such that the pullback $\pi^*\omega$ on $\pi^{-1}(X_\mathrm{reg})$ extends holomorphically to the whole of $Y$.

$\boldsymbol{\cdot}$ A desingularization $\pi:Y\to X$ of a symplectic variety $X$ is called a {\bf symplectic resolution} if the extended 2-form on $Y$ is still non-degenerate.
\end{dfn}

The only 2-dimensional symplectic singularities are ADE-singularities \cite[2.1]{Bea} and their symplectic resolutions are nothing but the minimal resolutions.

Symplectic varieties admit the following structure theorem by Kaledin.

\begin{thm}(\cite[Thm. 2.3]{K})\label{thm:Kaledin}
Let $X$ be a symplectic variety and $X\supset X_1\supset X_2\supset\cdots$ be the stratification obtained by successively taking the singular loci. Then each irreducible component $S$ (called a symplectic leaf) of $X_i\setminus X_{i+1}$ is a symplectic manifold and in particular even-dimensional. Moreover, at each point $x$ of a symplectic leaf $S$ there is a decomposition
$$\hat{X}_x\cong \mathcal{T}\times \hat{S}_x$$
where $\hat{X}_x$ and $\hat{S}_x$ are the formal completions of $X$ and $S$ at $x$ respectively and $\mathcal{T}$ is a formal symplectic variety (called a transversal slice of $S$). 
\end{thm}

When a symplectic leaf $S$ is of codimension 2 in $X$, the transversal slice $\mathcal{T}$ is 2-dimensional and thus it must have ADE-singularities. We will mainly treat 4-dimensional symplectic varieties in this paper. When we are given a symplectic resolution of a 4-dimensional symplectic variety, we will call the fiber of a 0-dimensional symplectic leaf a {\em central fiber}.

\begin{prop}\label{prop:sympA}
Assume $T=\mathbf{A}_n$ with $n\ge3$. Then $\mathrm{Hilb}^2(\Gamma_n)$ is a symplectic variety and has a unique 0-dimensional symplectic leaf $\{q\}$ in the sense of Theorem \ref{thm:Kaledin}.
\end{prop}

{\em Proof.} From the explicit calculation above, we know that $\mathrm{Hilb}^2(\Gamma_n)$ is a local complete intersection and its singular locus is of codimension 2. Therefore, $\mathrm{Hilb}^2(\Gamma_n)$ is normal by Serre's criterion. Moreover, the smooth part of $\mathrm{Hilb}^2(\Gamma_n)$ has a symplectic form outside a 2-dimensional locus since $\Gamma_n\setminus\{o\}$ has a symplectic form induced from $\mathbb{C}^2$ (cf. \cite{Fu}). Since $\mathrm{Hilb}^2(\Gamma_n)$ is birational to a symplectic manifold $\mathrm{Hilb}^2(\tilde{\Gamma}_n)$ where $\tilde{\Gamma}_n$ denotes the minimal resolution of $\Gamma_n$, the same argument as in the proof of \cite[Thm. 3.6]{Le} shows that the extension property of the symplectic form is satisfied. The second claim also follows from the calculation above.\qed

\vspace{3mm}

Next we show that successive blowing-ups give a symplectic resolution of $\mathrm{Hilb}^2(\Gamma_n)$ by explicit calculations. This can also be done by abstract argument (cf. Proposition \ref{prop:resolD}).

The singularities on the open subsets $U_i$ except $i=5$ are analytically locally just the product of $\mathbb{C}^2$ and ADE-singularities by Theorem \ref{thm:Kaledin}. Thus we concentrate on $U_5$ for $n\ge3$.

\begin{lem}
Let $\tau:\tilde{U}_5\to U_5$ be the blowing-up along the (reduced) singular locus. If $n=3$ or $4$, then $\mathrm{Sing}(\tilde{U}_5)$ is smooth. If $n\ge5$, then $\tilde{U}_5$ has a unique 0-dimensional symplectic leaf whose analytic germ is the same as that of $\mathrm{Hilb}^2(\Gamma_{n-2})$.

Moreover, for $n\ge3$, the exceptional locus of $\tau$ consists of two irreducible divisors, and $\tau^{-1}(q)$ is isomorphic to $\mathbb{P}^1\times\mathbb{P}^1$.
\end{lem}

{\em Proof.} Since we are interested in the singularity type at the origin, we may apply the following change of variables
$$x_1\mapsto x_1+x_2x_7',\;x_3\mapsto x_3+x_2x_9',\;x_5\mapsto x_5+x_2^2$$
so that the singular locus becomes $\{x_1=x_3=x_5=x_2^{n-1}-x_7'x_9'=0\}$. Let $\tilde{f}_3$ and $\tilde{f}_4$ be the polynomials obtained by applying the above transformation to $\tilde{f}_1$ and $\tilde{f}_2$ respectively. Put $t=x_2^{n-1}-x_7'x_9'$. Then one can check
{\footnotesize
$$\begin{aligned}\tilde{f}_4&=\sum_{k=1}^{\floor{n/2}}
2^{n-2k}\begin{pmatrix}n-k\\k\end{pmatrix}x_2^{n-2k}x_5^k-2^{n-1}(x_1x_9'+x_3x_7'-2x_2 t),\\
\tilde{f}_5:=\tilde{f}_3-x_2\tilde{f}_4&=\sum_{k=2}^{\floor{(n+1)/2}}2^{n-2k+1}
\begin{pmatrix}n-k\\k-1\end{pmatrix}x_2^{n-2k}x_5^k-2^{n-1}(x_1x_3+x_5t)
\end{aligned}$$}
by using the following formulas of binomial coefficients (see Remark \ref{rem1}.)
{\footnotesize
\begin{equation}\label{eq:2.1}
\sum_{i=k}^{\floor{(n+1)/2}}
\begin{pmatrix}n+1\\2i\end{pmatrix}
\begin{pmatrix}i\\k\end{pmatrix}
=2^{n-2k}\frac{n}{n-k}
\begin{pmatrix}n-k+1\\k\end{pmatrix},
\end{equation}

\begin{equation}\label{eq:2.2}
\sum_{i=k}^{\floor{n/2}}
\begin{pmatrix}n+1\\2i+1\end{pmatrix}
\begin{pmatrix}i\\k\end{pmatrix}
=2^{n-2k}\begin{pmatrix}n-k\\k\end{pmatrix}.
\end{equation}}

Let $\tilde{U}$ be the blowing-up of $U_5$ (with the coordinates transformed as above) along the singular locus $\{x_1=x_3=x_5=t=0\}$, and let $z_1,z_3,z_5,t'$ be the homogeneous coordinates of $\tilde{U}$ corresponding to $x_1,x_3,x_5,t$ respectively. Then similarly to the case of $\mathrm{Bl}_\Delta$, there exists only one affine chart $U:=\{z_5\ne0\}\subset\tilde{U}$ whose singular locus is singular. Set $x_i''=z_i/z_5,\;i=1,3$ and $t''=t'/z_5$. Then one sees that one can eliminate $t''$ and $U$ is defined by the following two equations in $\mathbb{C}^6=\mathrm{Spec}\,\mathbb{C}[x_1'',x_2,x_3'',x_5,x_7',x_9']$:
{\footnotesize
$$\begin{aligned}g_1&=\sum_{k=2}^{\floor{(n+1)/2}}
2^{n-2k+1}\begin{pmatrix}n-k-1\\k\end{pmatrix}x_2^{n-2k+1}x_5^{k-1}-2^{n-1}(x_1''x_3''x_5-x_2^{n-1}+x_7'x_9'),\\
g_2&=\sum_{k=1}^{\floor{n/2}}2^{n-2k}
\begin{pmatrix}n-k-1\\k-1\end{pmatrix}x_2^{n-2k}x_5^{k-1}-2^{n-1}(x_1''x_9'+x_3''x_7'-2x_1''x_2x_3'').
\end{aligned}$$}
Applying the transformations $x_7'\mapsto x_7'+2x_1''x_2$ and $x_9'\mapsto x_9'+2x_2x_3''$, the polynomials $g_1-2x_2g_2$ and $g_2$ become the same forms as $\tilde{f}_5$ and $\tilde{f}_4$ respectively with $n$ replaced by $n-2$. This proves the former part of the statement.

The latter part also follows from direct computation. Two irreducible exceptional divisors are defined as $\{x_5=x_7'-x_2x_1''=0\}$ and $\{x_5=x_9'-x_2x_3''=0\}$ in $U$. The central fiber in $\tilde{U}$ is defined by $z_1z_3-z_5t'=0$ in $\mathbb{P}^3=\mathrm{Proj}\,\mathbb{C}[z_1,z_3,z_5,t']$ and hence isomorphic to $\mathbb{P}^1\times\mathbb{P}^1$.
\qed

\vspace{3mm}

\begin{rem}\label{rem1}
The formulas (\ref{eq:2.1}) and (\ref{eq:2.2}) can be obtained, for example, using the Wilf-Zeilberger method (cf. \cite[Ch. 7]{PWZ}):\\
To adapt the notation to {\em ibid.}, we switch $i$ and $k$. For the formula (\ref{2.1}), it suffices to show that $\sum_{k=-\infty}^\infty F(n,k)=1$ for any $n$ where
{\footnotesize
$$F(n,k)=\begin{pmatrix}n+1\\2k\end{pmatrix}
\begin{pmatrix}k\\i\end{pmatrix}
\Bigg{/}\left\{2^{n-2i}\frac{n}{n-i}
\begin{pmatrix}n-i+1\\i\end{pmatrix}\right\}.$$}
This is done by finding rational functions $G(n,k)$ such that $\lim_{k\to \pm\infty}G(n,k)=0$ and
$$F(n+1,k)-F(n,k)=G(n,k+1)-G(n,k),\;\forall n,k.$$

In our case, if we set
{\footnotesize
$$R(n,k)=\frac{(i-k)(2k-1)}{(n-i)(n-2k+1)},$$}
then $G(n,k):=F(n,k)R(n,k)$ satisfies the above conditions.

For the formula (\ref{2.2}), the same method works by setting
{\footnotesize
$$R(n,k)=\frac{(i-k)(2k+1)}{(n-i)(n-2k)}. $$}

\vspace{-10mm}

\qed
\end{rem}

\vspace{3mm}

From the computation in the above proof, we also know that the general fiber of $\tau$ of the singular locus is a union of two copies of $\mathbb{P}^1$ intersecting transversally at one point if $n\ge2$ and is a single $\mathbb{P}^1$ if $n=1$. Moreover, the intersection of the central fiber (isomorphic to $\mathbb{P}^1\times\mathbb{P}^1$) and the singular locus is a union of 2 copies of $\mathbb{P}^1$ intersecting transversally at one point if $n\ge4$. If $n=3$, the intersection is a diagonal of $\mathbb{P}^1\times\mathbb{P}^1$.

\begin{prop}\label{prop:resolA}
$\floor{(n+1)/2}$ successive blowing-ups along singular loci give a symplectic resolution of $\mathrm{Hilb}^2(\Gamma_n)$. 
\end{prop}

{\em Proof.} Let $\psi:\mathcal{H}\to \mathrm{Hilb}^2(\Gamma_n)$ be the composite of the successive blowing-ups. Then $\mathcal{H}$ is smooth by the lemma above.

For the claim of symplecity, note that $\mathcal{H}$ is a canonical resolution of $\mathbf{A}_n$-singularities outside $q$. Therefore, $\mathcal{H}$ admits a natural symplectic form outside the fiber $F$ of $q$ and it extends to a 2-form $\omega$ on the whole of $\mathcal{H}$ since $\mathrm{Hilb}^2(\Gamma_n)$ is a symplectic variety (see Proposition \ref{prop:sympA}). Since $\mathcal{H}$ is smooth, the degeneracy locus of $\omega$ is a divisor obtained as the zero set of $\omega^{\wedge2}\in H^0(\mathcal{H},K_\mathcal{H})$. However, $F$ is of codimension 2 and hence $\omega$ is non-degenerate everywhere.
\qed

\vspace{3mm}

One of the main purposes of this article is to characterize the singularity of $\mathrm{Hilb}^2(\Gamma_n)$ by fibers of a symplectic resolution. Since direct calculation of the fiber $F$ of $q$ under $\psi:\mathcal{H}\to \mathrm{Hilb}^2(\Gamma_n)$ is hard, we use the results of Andreatta and Wi\'sniewski \cite{AW}. They describe the central fiber of the special symplectic resolution $\mathrm{Hilb}^2(\tilde{\Gamma}_n)$ of $\mathrm{Sym}^2(\Gamma_n)$ and then perform Mukai flops where $\tilde{\Gamma}_n$ is the minimal resolution of $\Gamma_n$.

Let us recall the results in \cite[\S6.4]{AW}. Let $\{C_1,\dots,C_n\}$ be the set of all  irreducible exceptional curves of $\tilde{\Gamma}_n\to \Gamma_n$ such that $C_i\cap C_{i+1}\ne\emptyset$ for $i=1,\dots,n-1$. Then the central fiber of the natural map $\pi_1:\mathrm{Sym}^2(\tilde{\Gamma}_n)\to\mathrm{Sym}^2(\Gamma_n)$ consists of $n(n+1)/2$ irreducible components: the images of $C_i\times C_j\subset(\tilde{\Gamma}_n)^{\times 2}\;(i\le j)$ in $\mathrm{Sym}^2(\tilde{\Gamma}_n)$. Let $P_{i,j}$ be the strict transforms of these components via the inverse of the Hilbert-Chow morphism $\pi_2:\mathrm{Hilb}^2(\tilde{\Gamma}_n)\to\mathrm{Sym}^2(\tilde{\Gamma}_n)$. Then $P_{i,i}$ is isomorphic to $\mathbb{P}^2$, $P_{i,i+1}$ is isomorphic to $\mathbb{P}^1\times\mathbb{P}^1$ blown up at one point, and another $P_{i,j}$ is  isomorphic to $\mathbb{P}^1\times\mathbb{P}^1$. The central fiber of $\pi=\pi_1\circ\pi_2$ consists of $P_{i,j}$'s and $n$ more irreducible components $Q_i\;(i=1,\dots,n)$, each of which is a $\mathbb{P}^1$-bundle over a conic in $P_{i,i}$. Also, $Q_i$ is isomorphic to $\Sigma_4$ and $Q_i\cap P_{i,i}$ is the unique negative section in $Q_i$ where $\Sigma_k=\mathbb{P}(\mathcal{O}_{\mathbb{P}^1}\oplus\mathcal{O}_{\mathbb{P}^1}(k))$ is the $k$-th Hirzebruch surface.

The successive blowing-ups $\psi:\mathcal{H}\to \mathrm{Hilb}^2(\Gamma_n)$ followed by the Hilbert-Chow morphism $\phi:\mathrm{Hilb}^2(\Gamma_n)\to\mathrm{Sym}^2(\Gamma_n)$ is a symplectic resolution of $\mathrm{Sym}^2(\Gamma_n)$ and thus it is obtained by a sequence of Mukai flops from $\mathrm{Hilb}^2(\tilde{\Gamma}_n)$ \cite[Thm. 1.2]{WW}.

\vspace{5mm}

\[\hspace{-25mm}
\xymatrix{
\mathcal{H} \ar[rr]^\psi \ar@{<.>}[d]_{\text{a sequence of Mukai flops}} &  & \mathrm{Hilb}^2(\Gamma_n) \ar[d]^{\phi} \\
\mathrm{Hilb}^2(\tilde{\Gamma}_n) \ar[r]^{\pi_2} & \mathrm{Sym}^2(\tilde{\Gamma}_n) \ar[r]^{\pi_1}  & \mathrm{Sym}^2(\Gamma_n)
}
\]

\vspace{7mm}

A Mukai flop is the following operation (see e.g. \cite[Ex. 2.6]{AW}). If one blows up a $\mathbb{P}^2$ in a 4-dimenional symplectic manifold $M$, then the exceptional divisor $E$ is isomorphic to the incidence variety
$$\{(p,l)\in\mathbb{P}^2\times(\mathbb{P}^2)^*\mid p\in l\}$$
where $(\mathbb{P}^2)^*$ is the projective space of lines in $\mathbb{P}^2$. Then one can blow down $E$ in the different direction (i.e., $E$ is projected to $(\mathbb{P}^2)^*$) to obtain a new symplectic manifold $M'$.  By symmetry, if one performs a Mukai flop again along $(\mathbb{P}^2)^*\subset M'$, one goes back to the original situation.

To determine how we should perform Mukai flops, we introduce the theory developed in \cite{AW} about the wall-and-chamber structure on the movable cone of a symplectic resolution. Let $E_0\subset \mathrm{Hilb}^2(\tilde{\Gamma}_n)$ be the irreducible $\pi_2$-exceptional divisor and set $\pi=\pi_2\circ\pi_1$. The other $\pi$-exceptional divisors $E_1,\dots,E_n$ are the strict transforms of the $\pi_1$-exceptional divisors which are the images of $C_i\times \tilde{\Gamma}_n\subset(\tilde{\Gamma}_n)^{\times 2}\,(i=1,\dots,n)$ in $\mathrm{Sym}^2(\tilde{\Gamma}_n)$. By abuse of notation, we also denote the strict transforms of $E_0,E_1,\dots,E_n$ in other symplectic resolutions (including $\mathcal{H}$) of $\mathrm{Sym}^2(\Gamma_n)$ by the same symbols. This makes sense since symplectic resolutions are isomorphic in codimension 1. Then $E_0,E_1,\dots,E_n$ form a basis of the relative Picard group $N^1(\pi)=[\mathrm{Pic}(\mathrm{Hilb}^2(\tilde{\Gamma}_n))/\pi^*\mathrm{Pic}(\mathrm{Sym}^2(\Gamma_n))]\otimes \mathbb{R}$ with real coefficients. Note that $N^1(\pi)$ and $N^1(\phi\circ\psi)$ are canonically isomorphic. 

Dually, the class $e_0,e_1,\dots,e_n$ of the general fibers of $\pi|_{E_0},\pi|_{E_1},\dots,\pi|_{E_n}$ respectively form a basis of the space $N_1(\pi)$ of curves which is the dual of $N^1(\pi)$ with respect to the intersection pairing. The intersection numbers $-(E_i.e_j)$ give (the direct sum of) Cartan matrices \cite[\S 4]{AW}.

Inside $N^1(\pi)$, the set of divisors $D$ such that the linear system $|mD|$ do not have fixed components for $m\gg0$ form a closed convex cone $\mathrm{Mov}(\pi)$ called the movable cone, and in fact $\mathrm{Mov}(\pi)$ is the dual cone of the cone generated by $e_0,e_1,\dots,e_n$ \cite[Thm. 4.1]{AW}. $\mathrm{Mov}(\pi)$ has a wall-and-chamber structure such that the set of chambers bijectively corresponds to the set of non-isomorphic symplectic resolutions. In our case the walls are defined by roots in the root system of type $\mathbf{A}_n$. More precisely, the walls are the hyperplanes that are orthogonal to the class $\lambda_{i,j}=e_0-(e_i+e_{i+1}+\cdots+e_j)\,(1\le i \le j\le n)$ \cite[Thm. 6.6]{AW}.

When we are given a symplectic resolution and a $\mathbb{P}^2$ inside it, we call the class $\lambda$ of a line in this $\mathbb{P}^2$ a {\em flopping class}. In this case the corresponding chamber has a facet defined as $\lambda^\perp$, and if we move to the adjacent chamber with respect to $\lambda^\perp$, then only this $\mathbb{P}^2$ is flopped and a line in the flopped $\mathbb{P}^2$ gives the class $-\lambda$. 

Now we can determine which chamber corresponds to $\mathcal{H}$. First note that the central fiber of $\phi\circ\psi$ contains a component $P$ which is the strict transform of the central fiber of $\phi$. Since the central fiber of $\phi$ (which is isomorphic to $\mathbb{P}^2$ as we calculated above) intersects with the singular locus of $\mathrm{Hilb}^2(\Gamma_n)$ along a curve, $P$ is also isomorphic to $\mathbb{P}^2$. Also, $P$ is not contained in $E_1,\dots,E_n$. Therefore, the flopping class of $P$ intersects with $E_1,\dots,E_n$ non-negatively. One can easily check that $\lambda_{1,n}$ is a unique such class. There is a unique chamber whose corresponding resolution has the flopping class $\lambda_{1,n}$. This chamber is defined as $\{D\in\mathrm{Mov}(\pi)\mid (D.\lambda_{1,n})\ge0\}$.

In \cite{AW} the authors actually perform Mukai flops. They describe the case $n=6$, and the general cases are similar. The above-mentioned chamber corresponds to the final diagram of \cite[\S6]{AW}. To obtain $\mathcal{H}$ from $\mathrm{Hilb}^2(\tilde{\Gamma}_n)$, one should perform Mukai flops to $P_{i,j}$'s, for example, in the following order
$$P_{1,1}\to P_{2,2}\to P_{2,1}\to P_{3,3}\to P_{3,2}\to P_{3,1}\to P_{4,4}\to\cdots \to P_{n,2}\to P_{n,1}$$
where, by abuse of notation, $P_{i,j}$ denotes its strict transforms under Mukai flops. We also denote its strict transforms by $Q_i$. If $P_{i,j}$ is isomorphic to $\mathbb{P}^2$ and is flopped, the same symbol also denotes the dual space $(\mathbb{P}^2)^*$ in the flopped resolution.

\begin{figure}[H]
\centering
\includegraphics[
scale=0.3
]{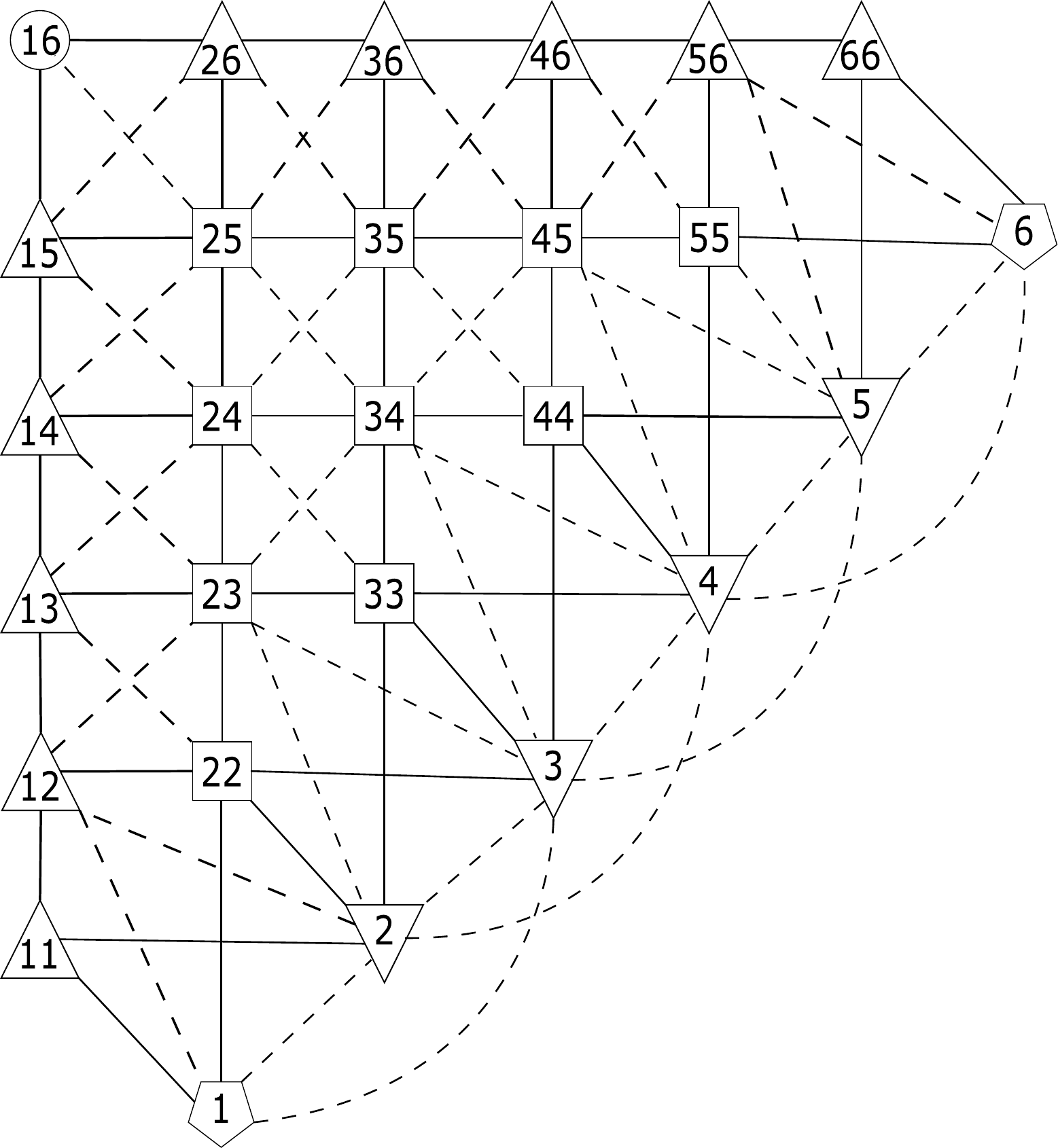}
\caption{The central fiber of $\phi\circ\psi$ for $\mathbf{A}_6$}
\label{fig1}
\end{figure}

In the diagrams of \cite[\S6]{AW}, however, (the strict transforms of) $Q_i$'s are omitted for the reason that these components will never become $\mathbb{P}^2$ under any sequence of Mukai flops and thus will never be flopped. One can check that, under the above sequence of Mukai flops, $Q_1$ and $Q_n$ become $\Sigma_3$'s, and $Q_2,\dots,Q_{n-1}$ become $\Sigma_2$'s in the final diagram (see Figure \ref{fig1} for $n=6$).

In Figure \ref{fig1}, the symbols with the numbers $ij$ (resp. $i$) stand for the components $P_{i,j}$ (resp. $Q_i$). The shapes of the symbols stand for the isomorphism classes of the components: circles, (upright) squares, upright triangles, inverted triangles, and inverted pentagons indicate $\mathbb{P}^2,\,\mathbb{P}^1\times\mathbb{P}^1,\, \Sigma_1,\,\Sigma_2$ and $\Sigma_3$ respectively. If two components intersect along a curve, which is always isomorphic to $\mathbb{P}^1$, the corresponding symbols are joined by a solid line. Similarly, a dotted line corresponds to a 0-dimensional intersection, which is always a single point.

Since the central fiber of $\phi$, whose strict transform in $\mathcal{H}$ is $P_{1,n-1}$, and $\mathrm{Sing}(\mathrm{Hilb}^2(\Gamma_n))$ intersects along two projective lines meeting at one point, we see that the central fiber $F$ of $\psi$ corresponds to the subdiagram consisting of all the squares and the inverted triangles. Thus the fiber consists of $n(n-1)/2$ copies of $\mathbb{P}^1\times\mathbb{P}^1$ and $n-2$ copies of $\Sigma_2$.

\vspace{3mm}

\begin{rem}\label{rem2}
The central fiber $F$ does not contain $\mathbb{P}^2$, and we can say that $\mathrm{Hilb}^2(\Gamma_n)$ has a unique symplectic resolution. In general uniqueness of symplectic resolutions of a 4-dimensional symplectic variety is equivalent to nonexistence of $\mathbb{P}^2$ in the central fiber of some symplectic resolution since any two symplectic resolutions are related by Mukai flops \cite{WW} and conversely every $\mathbb{P}^2$ can be flopped to give a different symplectic resolution.
\qed
\end{rem}

We can determine which components in $F$ are contained in which exceptional divisors. When $P_{i,j}$ is flopped, its flopping class is $-\lambda_{i,j}$ and turns into $\lambda_{i,j}$. Thus calculating the intersection number $(\lambda_{i,j}. E_k)$ shows whether or not $P_{i,j}$ is contained in $E_k$ after the flop. In the final diagram for the case $\mathbf{A}_n$ (see Figure \ref{fig1}), a component $P_{i,j}$ in $F$ is contained in $E_k\,(k=1,\dots,n)$ if and only if $i=k+1$ or $j=k-1$. Also, $Q_i$ in $F$ is contained in $E_k\,(k=1,\dots,n)$ if and only if $i=k$.

The configuration of the components of $F$ is explained as follows. The intersection $F_k:=F\cap E_k\,(k=1,\dots,n)$ is a $\mathbb{P}^1$-bundle over the Dynkin tree of $\mathbb{P}^1$'s of type $\mathbf{A}_{n-2}$. Moreover, $Q_i\,(i=2,\dots,n-1)$ intersects with $P_{i,i}$ along the $(-2)$-curve, which is a diagonal (i.e. an irreducible curve $C$ with $C^2=2$) in $P_{i,i}\cong\mathbb{P}^1\times\mathbb{P}^1$. Every 0-dimensional intersection is realized as the intersection of two 1-dimensional intersections. Thus one can recover $F$ up to isomorphism from these data.

\subsection{Case $T=\mathbf{D}_n$}\label{2.2}

In this subsection we consider the case when $T=\mathbf{D}_n$. Then the surface $\Gamma_n:=\Gamma$ is the quotient of $\mathbb{C}^2$ by the action of the binary dihedral group $BD_{n-2}$ of order $4(n-4)$ generated by {\footnotesize $\begin{pmatrix}\zeta&0\\0&\zeta^{-1}\end{pmatrix}$} and {\footnotesize $\begin{pmatrix}0&1\\-1&0\end{pmatrix}$} where $\zeta\in\mathbb{C}$ is a primitive $2(n-2)$-th root of unity. Therefore, $\mathrm{Sym}^2(\Gamma_n)$ is the quotient of $\mathrm{Spec}\,(\mathbb{C}[a_1,a_2]\otimes\mathbb{C}[b_1,b_2])=\mathbb{C}^2\times\mathbb{C}^2$ by the natural action of the wreath product $G=BD_{n-2}\wr \mathfrak{S}_2$.
The ring of invariants by the normal subgroup $(BD_{n-2})^{\times 2}\lhd G$ is generated by
{\footnotesize
$$u_1=a_1^2a_2^2,\;v_1=a_1^{2n-4}+a_2^{2n-4},\;w_1=a_1^{2n-3}a_2-a_1a_2^{2n-3},$$
$$u_2=b_1^2b_2^2,\;v_2=b_1^{2n-4}+b_2^{2n-4},\;w_2=b_1^{2n-3}b_2-b_1b_2^{2n-3}.$$}
These generators satisfy $u_i(v_i^2-4u_i^{n-2})-w_i^2=0,\,i=1,2$. The ring of invariants by $G$ is generated by the following elements:
{\footnotesize
$$\begin{aligned}x_1=u_1+u_2,\;
x_2&=v_1+v_2,\;
x_3=w_1+w_2,\\
x_4=(u_1-u_2)^2,\;
x_5&=(v_1-v_2)^2,\;
x_6=(w_1-w_2)^2,\\
x_7=(u_1-u_2)(v_1-v_2),\;
x_8&=(u_1-u_2)(w_1-w_2),\;
x_9=(v_1-v_2)(w_1-w_2).
\end{aligned}$$}

\begin{prop}\label{prop:relation2}
The ideal of relations of $x_i$'s is generated by the following ten elements of $\mathbb{C}[x_1,x_2,\dots,x_9]$:
{\footnotesize
$$\begin{aligned}
f_1=&\sum_{i=0}^{\floor{(n-1)/2}}
\begin{pmatrix}n-1\\2i\end{pmatrix}x_1^{n-2i-1}x_4^i+2^{n-6}
(-x_1x_2^2+2x_3^2-x_1x_5+2x_6-2x_2x_7),\\
f_2=&\sum_{i=0}^{\floor{n/2}-1}
\begin{pmatrix}n-1\\2i+1\end{pmatrix}x_1^{n-2i-2}x_4^{i+1}+2^{n-6}
(-x_2^2x_4-x_4x_5-2x_1x_2x_7+4x_3x_8),\\
f_3=&x_4x_5-x_7^2,\;f_4=x_7x_8-x_4x_9,\;f_5=x_4x_6-x_8^2,\\
f_6=&\sum_{i=0}^{\floor{n/2}-1}
\begin{pmatrix}n-1\\2i+1\end{pmatrix}x_1^{n-2i-2}x_4^i x_7+2^{n-6}
(-2x_1x_2x_5-x_2^2x_7-x_5x_7+4x_3x_9),\\
\end{aligned}$$
$$\begin{aligned}
f_7=&x_5x_8-x_7x_9,\\
f_8=&\sum_{i=0}^{\floor{n/2}-1}
\begin{pmatrix}n-1\\2i+1\end{pmatrix}x_1^{n-2i-2}x_4^i x_8+2^{n-6}
(4x_3x_6-x_2^2x_8-2x_1x_2x_9-x_7x_9),\\
f_9=&x_6x_7-x_8x_9,\;f_{10}=x_5x_6-x_9^2.
\end{aligned}$$}
\end{prop}

{\em Proof.} We can apply the same argument as in the proof of Proposition \ref{prop:relation1}.
\qed

\vspace{3mm}

We calculate the blowing-up of $\mathrm{Sym}^2(\Gamma_n)$ along the diagonal in a similar fashion to the case $T=\mathbf{A}_n$. We use the same notations as in Subsection \ref{2.1}. 

Similarly to the case $T=\mathbf{A}_n$, the central fiber of the blowing-up is isomorphic to $\mathbb{P}^2$ and the singular locus $\mathrm{Sing}(\mathrm{Bl}_\Delta)$ is isomorphic to the blowing-up $\mathrm{Bl}_o (\Gamma_n)$ of $\Gamma_n$ at the origin. In particular it has one $\mathbf{A}_1$-singularity and one $\mathbf{D}_{n-2}$-singularity if $n\ge5$, and has three $\mathbf{A}_1$-singularities if $n=4$ . For all $n\ge4$, the affine chart $U_5$ is defined by the following two equations in $\mathbb{C}^6=\mathrm{Spec}\,\mathbb{C}[x_1,x_2,x_3,x_5,x_7',x_9']$:
{\footnotesize
$$\begin{aligned}\tilde{f}_1=&\sum_{i=0}^{\floor{(n-1)/2}}
\begin{pmatrix}n-1\\2i\end{pmatrix}x_1^{n-2i-1}x_5^ix_7'^{2i}+2^{n-6}
(-x_1x_2^2+2x_3^2-x_1x_5+2x_5x_9'^2-2x_2x_5x_7'),\\
\tilde{f}_2=&\sum_{i=0}^{\floor{n/2}-1}
\begin{pmatrix}n-1\\2i+1\end{pmatrix}x_1^{n-2i-2}x_5^ix_7'^{2i+1}+2^{n-6}
(-2x_1x_2-x_2^2x_7'-x_5x_7'+4x_3x_9').
\end{aligned}$$}
The origin $q_1\in U_5$ corresponds to (one of) the $\mathbf{A}_1$-singular point(s). For $n=4$, the analytic germs of the other two $\mathbf{A}_1$-singular points in $\mathrm{Bl}_\Delta$ is the same as that of $q_1$ since the three singular points in $\mathrm{Bl}_o (\Gamma_n)$ are symmetric. We assume $n\ge5$. By abuse of notation, we set $x_i'=y_i/y_4,\;i=4,\dots,9$. Then the affine chart $U_4$ is defined by the following two equations in $\mathbb{C}^6=\mathrm{Spec}\,\mathbb{C}[x_1,x_2,x_3,x_4,x_7',x_8']$:
{\footnotesize
$$\begin{aligned}\tilde{f}_1=&\sum_{i=0}^{\floor{(n-1)/2}}
\begin{pmatrix}n-1\\2i\end{pmatrix}x_1^{n-2i-1}x_4^i+2^{n-6}
(-x_1x_2^2+2x_3^2-x_1x_4x_5+2x_4x_8'^2-2x_2x_4x_7'),\\
\tilde{f}_2=&\sum_{i=0}^{\floor{n/2}-1}
\begin{pmatrix}n-1\\2i+1\end{pmatrix}x_1^{n-2i-2}x_4^i+2^{n-6}
(-2x_1x_2x_7'-x_2^2-x_4x_7'^2+4x_3x_8').
\end{aligned}$$}
The origin $q_2\in U_4$ corresponds to the $\mathbf{D}_{n-2}$-singular point. We can prove the following in the same way as Proposition \ref{prop:sympA}.

\begin{prop}\label{prop:sympD}
Assume $T=\mathbf{D}_n$ with $n\ge4$. Then $\mathrm{Hilb}^2(\Gamma_n)$ is a symplectic variety. The 0-dimensional symplectic leaf in the sense of Theorem \ref{thm:Kaledin}. consists of three points (resp. the two points $q_1$ and $q_2$) if $n=4$ (resp. $n\ge5$). 
\end{prop}

\vspace{3mm}

Next we consider a symplectic resolution $\psi:\mathcal{H}\to\mathrm{Hilb}^2(\Gamma_n)$ for $T=\mathbf{D}_n$. As will see later in Proposition \ref{prop:resolD}, in fact the resolution $\psi$ uniquely exists and it is obtained by successive blowing-ups of singular loci as well as the $\mathbf{A}_n$-case. We admit the existence of $\psi$ for the moment. Then $\mathcal{H}$ is also a symplectic resolution of $\mathrm{Sym}^2(\Gamma_n)$ and thus it is obtained by performing Mukai flops to $\mathrm{Hilb}^2(\tilde{\Gamma}_n)$ where $\tilde{\Gamma}_n$ is the minimal resolution of $\Gamma_n$.

\begin{figure}[H]
\centering
\includegraphics[scale=0.3
]{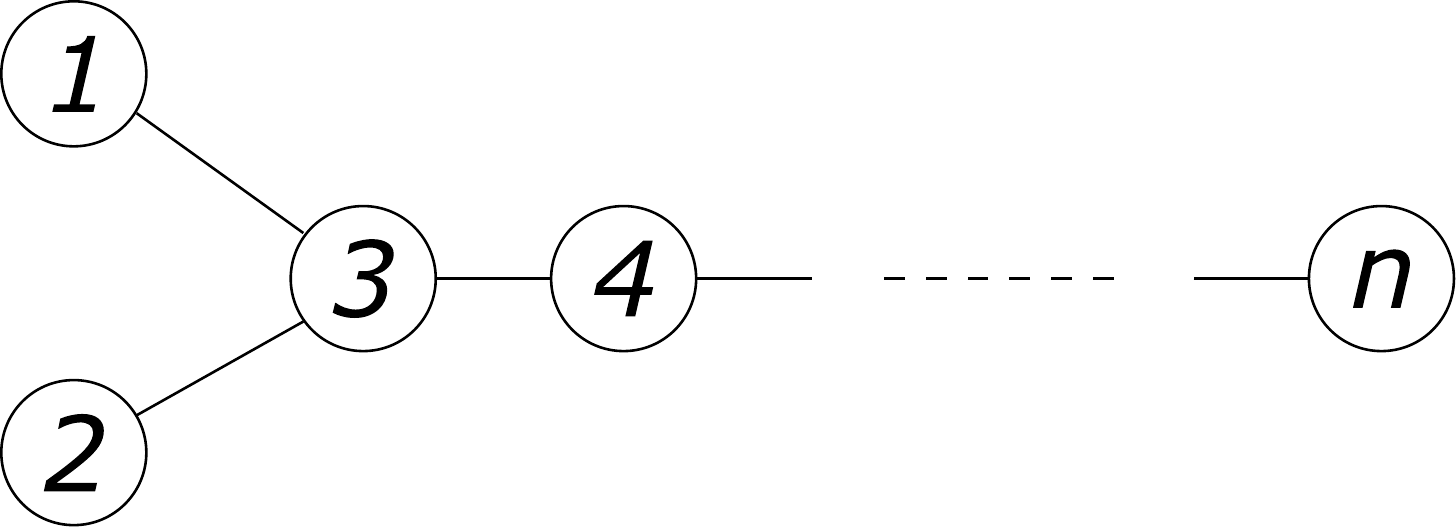}


\caption{The dual graph of the exceptional curves of type $\mathbf{D}_n$}
\label{fig2}
\end{figure}

To describe the Mukai flops, let $\{C_1,\dots,C_n\}$ be the set of the exceptional curves. Here each $C_i$ is chosen so that $C_i$ corresponds to the $i$-th vertex in the dual graph in Figure \ref{fig2}.

We define the surfaces $P_{i,j}$ and $Q_i$ in $\mathrm{Hilb}^2(\tilde{\Gamma}_n)$ in the same way as the cases $T=\mathbf{A}_n$. Similarly to $\mathbf{A}_n$, we can consider the movable cone and the chamber structure on it. We should just replace the root system $\mathbf{A}_n$ by $\mathbf{D}_n$. In this case we regard $e_i$ as the simple root that corresponds to the $i$-th vertex of the Dynkin diagram of $\mathbf{D}_n$ in Figure \ref{fig2}. Although the result about the chamber structure is shown only for $\mathbf{A}_n$ in \cite{AW}, the analogues for $\mathbf{D}_n$ and $\mathbf{E}_n$ (and even for higher dimensional cases) are shown by Bellamy \cite[4.2]{Bel}. In any case the chamber corresponding to $\mathcal{H}$ is described by the highest root $\alpha_h$, that is, it is defined by the hyperplane $(e_0-\alpha_h)^\perp$ in the movable cone (see the proof of Proposition \ref{prop:resolD}). To reach this chamber from $\mathrm{Hilb}^2(\tilde{\Gamma}_n)$, one must cross the same number of walls as the number of the positive roots. In $\mathbf{D}_n$ case, we need $n(n-1)$ flops.

We can write a positive root for $\mathbf{D}_n$ as $\alpha=\sum_{i=1}^n c_i e_i$ with $c_i\in\{0,1,2\}$. Then we define a total order $<$ on the set of the positive roots as follows. Two positive roots $\alpha=\sum_{i=1}^n c_i$ and $\alpha'=\sum_{i=1}^n c'_i$ satisfy $\alpha<\alpha'$ if and only if there exists $j\in\{1,\dots,n\}$ such that $c_j<c'_j$ and $c_k=c'_k$ for all $k>j$:
$$\alpha_1<\alpha_2<\alpha_3<\alpha_1+\alpha_3<\alpha_2+\alpha_3<\alpha_1+\alpha_2+\alpha_3<\alpha_4<\alpha_3+\alpha_4<\cdots<\alpha_h.$$
Then we can perform Mukai flops so that we cross the hyperplanes $(e_0-\alpha)^\perp$ in this order.

Through the crossings of the hyperplanes, the components of the central fiber are flopped in the following order:
{\footnotesize
$$\begin{aligned}
&P_{1,1}\to P_{2,2}\to P_{3,3}\to P_{1,3}\to P_{2,3}\to P_{1,2}\to P_{4,4}\to P_{3,4}\to P_{1,4}\to P_{2,4}\to P_{3,3}(2)\\
\to &Q_3\,(\mathbf{D}_4 \text{ ends here})\to P_{5,5}\to P_{4,5}\to P_{3,5}\to P_{1,5}\to P_{2,5}\to P_{3,4}(2)\to P_{4,4}(2)\to Q_4\,(\mathbf{D}_5 \text{ ends here})\\
\to &P_{6,6}\to P_{5,6}\to P_{4,6}\to P_{3,6}\to P_{1,6}\to P_{2,6}\to P_{3,5}(2)\to P_{4,5}(2)\to P_{5,5}(2)\to Q_5\,(\mathbf{D}_6 \text{ ends here})\\
\to &P_{7,7}\to P_{6,7}\to P_{5,7}\to P_{4,7}\to P_{3,7}\to P_{1,7}\to P_{2,7}\to P_{3,6}(2)\to P_{4,6}(2)\to P_{5,6}(2)\to P_{6,6}(2)\\
\to &Q_6\,(\mathbf{D}_7 \text{ ends here})\to P_{8,8}\to P_{7,8}\to\cdots\to P_{n-1,n-1}(2)\to Q_{n-1}
\end{aligned}$$}
where a number $k$ in parentheses means the $k$-th time flop of the component. Unlike $\mathbf{A}_n$, one component can be flopped more than once. Also, $Q_i$'s can be flopped.

\begin{figure}[H]
\centering
\includegraphics[scale=0.25
]{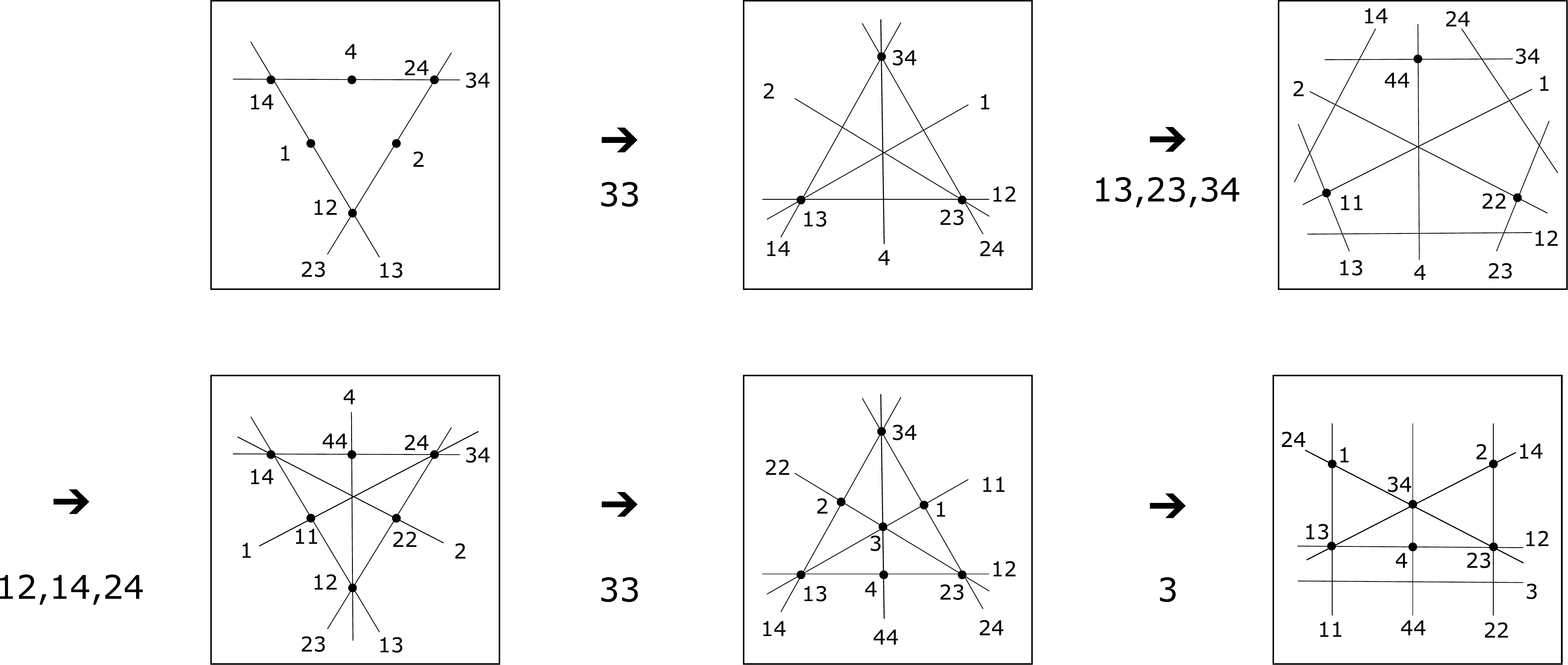}
\caption{Intersections of $P_{3,3}$ with other components for $\mathbf{D}_4$}
\label{fig3}
\end{figure}

The flops for $\mathbf{D}_n$ are more involved than the case $T=\mathbf{A}_n$. We show how the flops are done for the component $P_{3,3}$ for $\mathbf{D}_4$ as an example (Figure \ref{fig3}). For the other components, description of the flops is easier.

In Figure \ref{fig3}, all the intersections of $P_{3,3}$ with other components (except $Q_3$ in the first four diagrams, for readability) are described. A line or a point which is labeled with a number $ij$ (resp. $i$) in each diagram indicates the intersection with $P_{i,j}$ (resp. $Q_i$). A line (resp. a point) corresponds to a 1-dimensional (resp. 0-dimensional) intersection. In fact any 1-dimensional intersection is isomorphic to $\mathbb{P}^1$. The numbers under the arrows mean that the components labeled with them are flopped. For convenience, the order of the flops is a little different from the one that we gave above.

In the beginning, $P_{3,3}$ is isomorphic to $\mathbb{P}^2$. The intersection with $Q_3$ (which is omitted in the first diagram) is the conic passing through the points $1,2$ and $4$. $P_{3,3}$ itself is flopped and we move to the second diagram. Note that the roles of lines and points are switched between the first and the second diagrams. The intersection with $Q_3$ in the second diagram is the conic passing through the points $13,23$ and $34$. The process from the second to the fourth diagrams is the standard Cremona transformation, and thus $P_{3,3}$ again becomes $\mathbb{P}^2$. Also, the intersection with $Q_3$ becomes the line passing through the points $11,22$ and $44$. Then $P_{3,3}$ is itself flopped again and finally becomes a $\mathbb{P}^2$ blown up at one point at the sixth diagram.

\begin{figure}[H]
\centering
\includegraphics[scale=0.3
]{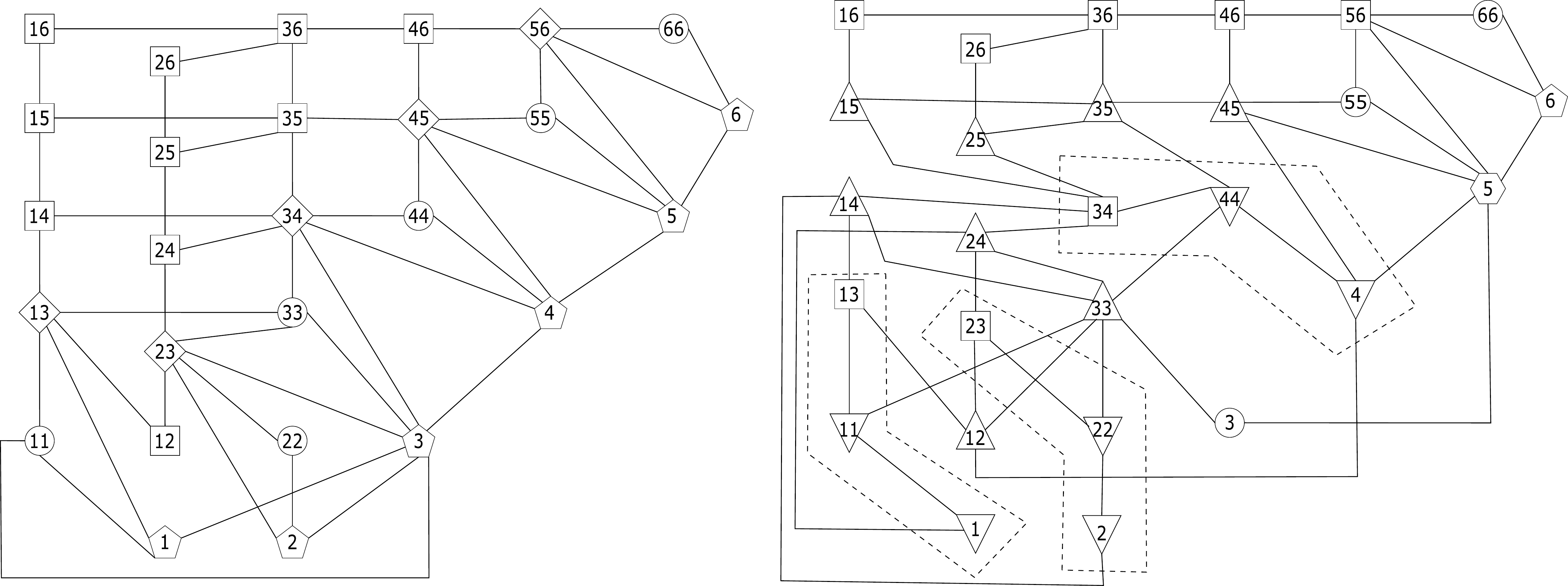}
\caption{The initial diagram and the central fibers of $\psi$ for $\mathbf{D}_4$}
\label{fig4}
\end{figure}

\begin{figure}[H]
\centering
\includegraphics[scale=0.3
]{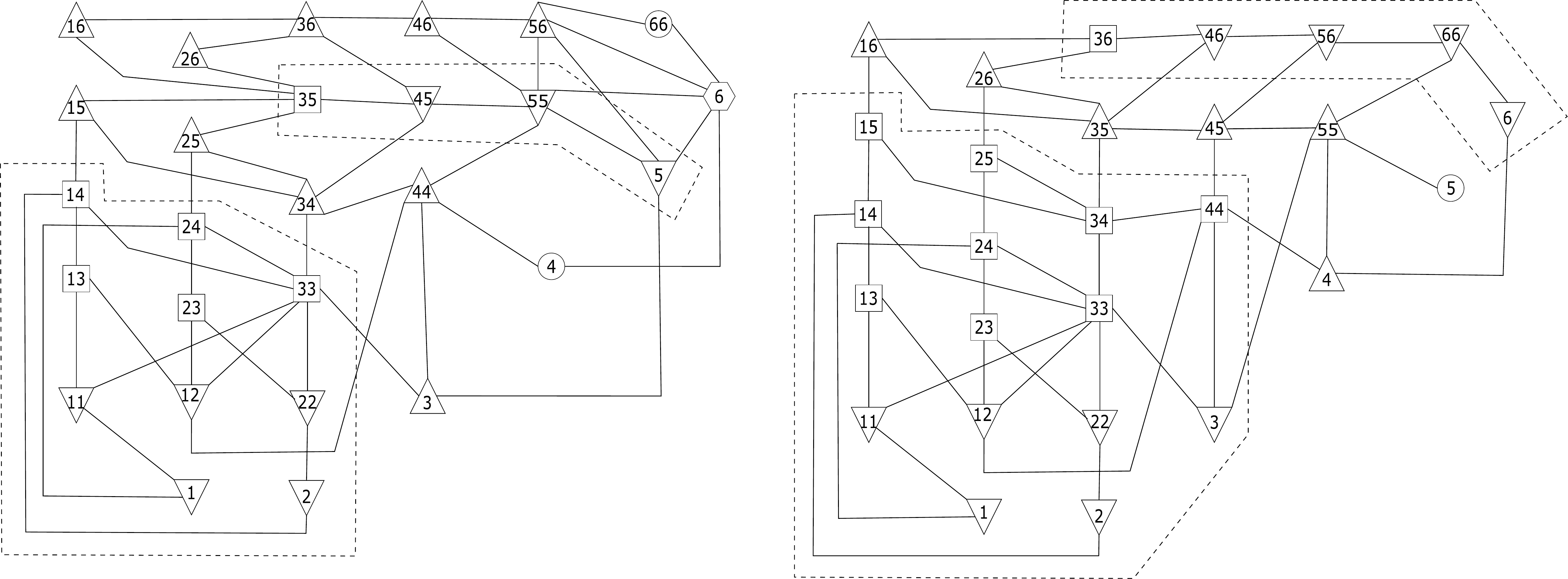}
\caption{The central fibers of $\psi$ for $\mathbf{D}_5$ and $\mathbf{D}_6$}
\label{fig5}
\end{figure}

We give the initial and the final diagrams of the flopping process below for $n=4,5,6$ and indicate the central fiber of $\pi$ (Figure \ref{fig4}, \ref{fig5}). The presentation is basically same as Figure \ref{fig1}. However, we omit 0-dimensional intersections between components in order to improve readability. Also, new symbols are there: upright pentagons (resp. rotated squares) stand for $\Sigma_4$ (resp. $\mathbb{P}^2$ blown up at two points). A hexagon indicates some smooth rational surface which is not isomorphic to the surfaces of the other symbols. The subdiagrams that corresponds to the central fibers of $\psi$ are enclosed by dotted lines.


Similarly to $\mathbf{A}_n$, we can determine the intersection of the central fibers of $\psi$ and the exceptional divisors. Assume $n\ge5$ and set $F_i=\psi^{-1}(q_i),\,i=1,2$ (see Proposition \ref{prop:sympD}). Then we have
{\footnotesize
$$F_1\cap E_k=
\begin{cases}P_{3,n}&\text{if }k=1,2\\
P_{k+1,n}&\text{if }k=3,\dots,n-1\\
Q_n&\text{if }k=n
\end{cases}$$}
and
{\footnotesize
$$F_2\cap E_k=
\begin{cases}P_{2,3}\cup P_{2,4}\cup P_{2,5}\cup\cdots\cup P_{2,n-1}\cup Q_1&\text{if }k=1\\
P_{1,3}\cup P_{1,4}\cup P_{1,5}\cup\cdots\cup P_{1,n-1}\cup Q_2&\text{if }k=2\\
P_{1,1}\cup P_{2,2}\cup P_{3,3}\cup P_{3,4}\cup\cdots\cup P_{3,n-2}&\text{if }k=3\\
P_{1,2}\cup P_{1,3}\cup P_{2,3}\cup P_{4,4}\cup P_{4,5}\cup\cdots\cup P_{4,n-2}&\text{if }k=4\\
\begin{aligned}
&P_{1,k-1}\cup P_{2,k-1}\cup P_{3,k-2}\cup P_{4,k-2}\cup \cdots\cup P_{k-2,k-2}\\
&\cup P_{k,k}\cup P_{k,k+1}\cup \cdots\cup P_{k,n-2}\cup Q_{k-2}
\end{aligned}
&\text{if }5\le k<n\\
P_{1,n-1}\cup P_{2,n-1}\cup P_{3,n-2}\cup P_{4,n-2}\cup \cdots\cup P_{n-2,n-2}&\text{if } k=n.
\end{cases}$$}

In any case, $F_2\cap E_k$ is a $\mathbb{P}^1$-bundle over a Dynkin tree of $\mathbb{P}^1$'s of type $\mathbf{D}_{n-2}$.

For $n=4$, let $F$ be the (disjoint) union of the three central fibers. Then we have
{\footnotesize
$$F\cap E_1=P_{2,3}\sqcup P_{3,4}\sqcup Q_1,\;
F\cap E_2=P_{1,3}\sqcup P_{3,4}\sqcup Q_2$$
$$F\cap E_3=P_{1,1}\sqcup P_{2,2}\sqcup P_{4,4},\;
F\cap E_4=P_{1,3}\sqcup P_{2,3}\sqcup Q_4.$$}

The configuration of the components in $F_i$ or $F$ is analogous to the $\mathbf{A}_n$-case. A component of $F$ which is isomorphic to $\mathbb{P}^1\times\mathbb{P}^1$ intersects with another component along a ruling or a diagonal (cf. Subsection \ref{2.1}). Unlike $\mathbf{A}_n$, two components which are isomorphic to $\Sigma_2$ can have a 1-dimensional intersection. This is the $(-2)$-curve of one component and is a section $C$ with $C^2=2$ of the other component. These data are sufficient to recover $F_i$ and $F$ up to isomorphism.

Using the above description of $F_1$ and $F_2$ and Theorem \ref{thm:main}, which we will prove in Section \ref{3}, we can show that blowing-ups along the singular loci give a symplectic resolution of $\mathrm{Hilb}^2(\Gamma_n)$. (We directly calculated the blowing-ups for $T=\mathbf{A}_n$ in Subsection \ref{2.1}, but calculation for $T=\mathbf{D}_n$ is harder and we adopt another approach.)

\begin{prop}\label{prop:resolD}
There exists a unique symplectic resolution of $\mathrm{Hilb}^2(\Gamma_n)$, and it is obtained by successively blowing up the singular loci of $\mathrm{Hilb}^2(\Gamma_n)$.
\end{prop}

{\em Proof.} Since $\mathrm{Hilb}^2(\tilde{\Gamma}_n)\to \mathrm{Sym}^2(\Gamma_n)$ is a symplectic resolution, $\mathrm{Hilb}^2(\tilde{\Gamma}_n)$ is a Mori dream space relative to $\mathrm{Sym}^2(\Gamma_n)$ (cf. \cite[Thm. 3.2]{AW}). The central fiber $P\cong\mathbb{P}^2$ of $\mathrm{Hilb}^2(\Gamma_n)\to \mathrm{Sym}^2(\Gamma_n)$ is two-dimensional, and thus the rational map $\mathrm{Hilb}^2(\tilde{\Gamma}_n)\dashrightarrow \mathrm{Hilb}^2(\Gamma_n)$ becomes a morphism $\psi:\mathcal{H}\to \mathrm{Hilb}^2(\Gamma_n)$ after some flops by the property of a Mori dream space \cite[1.11]{HK}. Since Mukai flops preserve smoothness and symplecity, $\mathcal{H}$ is a symplectic resolution of $\mathrm{Hilb}^2(\Gamma_n)$. This resolution corresponds to the chamber in the movable cone of $\mathrm{Hilb}^2(\tilde{\Gamma}_n)$ defined by the highest root since $P$ is not contained in $E_1,\dots,E_n$. In particular, $\psi$ is unique.

Let $\mathcal{H}_1$ be the blowing-up of $\mathrm{Hilb}^2(\Gamma_n)$ along the singular locus. Similarly as above, the central fibers are two-dimensional, and $\psi$ factors through $\mathcal{H}_1$ after some flops. Also, we see that $\mathcal{H}_1$ is normal since $\mathrm{Hilb}^2(\Gamma_n)$ is normal. (From this we know that $\mathcal{H}_1$ is also a possibly non-$\mathbb{Q}$-factorial Mori dream space, cf. \cite[\S 9]{Ok}.) Note that $F_1$ and $F_2$ contain no $\mathbb{P}^2$ and hence $\psi$ admits no flops (cf. Remark \ref{rem2}). Thus, the rational map $\psi_1:\mathcal{H}\dashrightarrow \mathcal{H}_1$ is a morphism, and $\mathcal{H}_1$ is a symplectic variety.

Considering the shapes of fibers of $\psi_1$ using the above description of $F_1$ and $F_2$, we see that the central fibers of $\psi_1$, if they exist, are isomorphic to the central fibers of a symplectic resolution of $\mathrm{Hilb}^2(\Gamma')$ with $\Gamma'$ of type $\mathbf{D}_{n-2}$. (Note that the exceptional divisor of $\mathcal{H}_1\to \mathrm{Hilb}^2(\Gamma_n)$ corresponds to $E_{n-1}$.) By Theorem \ref{thm:main}, $\mathcal{H}_1$ has the same singularities as $\mathrm{Hilb}^2(\Gamma')$ if it has 0-dimensional symplectic leaves. Therefore, by induction on relative Picard numbers, it follows that repeating blowing-ups gives the symplectic resolution.
\qed

\subsection{Case $T=\mathbf{E}_n$}\label{2.3}

By using software such as SINGULAR \cite{GPS}, we can calculate the defining equation of $\mathrm{Hilb}^2(\Gamma_n)$ for $T=\mathbf{E}_n$. Then one sees that the singular locus of $\mathrm{Hilb}^2(\Gamma_n)$ is isomorphic to the blowing-up $\mathrm{Bl}_o(\Gamma_n)$ as well as the cases $T=\mathbf{A}_n,\mathbf{D}_n$. The equation of  $\mathrm{Hilb}^2(\Gamma_n)$ in some affine open subset in $\mathrm{Spec}\,\mathbb{C}[x_1,x_2,x_3,x_4,x_5,x_6]$ whose origin is the unique 0-dimensional symplectic leaf is given as follows:
{\footnotesize
$$\begin{aligned}
&n=6\hspace{5mm}
\begin{cases}
&x_4x_5^3+ 2x_1^3+ 3x_2^2x_5+2x_1x_4+4x_3x_6,\\
&x_1^4+6x_2x_4x_5^2+2x_2^3+6x_1^2x_4+4x_4x_6^2+ 4x_3^2+x_4^2. \end{cases}\\
&n=7\hspace{5mm}
\begin{cases}
&x_1^3x_5+2x_4x_5^3+3x_1^2x_2+6x_2^2x_5+3x_1x_4x_5+x_2x_4+8x_3x_6,\\
&x_1^3x2+3x_1^2x_4x_5+6x_2x_4x_5^2+2x_2^3+3x_1x_2x_4+x_4^2x_5
+4x_4x_6^2+4x_3^2.
\end{cases}\\
&n=8\hspace{5mm}
\begin{cases}
&5x_1^4+4x_4x_5^3+10x_1^2x_4+12x_2^2x_5+x_4^2+16x_3x_6,\\ &x_1x_4x_5^3-10x_1^3x_4+3x_1x_2^2x_5-15x_2x_4x_5^2-5x_2^3
-6x_1x_4^2+ 4x_1x_3x_6-10x_4x_6^2-10x_3^2.
\end{cases}\\
\end{aligned}$$}

Note that in any case $\mathrm{Hilb}^2(\Gamma_n)$ is a local complete intersection. Therefore, we can prove the following similarly to the previous cases.

\begin{prop}\label{prop:sympE}
Assume $T=\mathbf{E}_n$ with $n=6,7,8$. Then $\mathrm{Hilb}^2(\Gamma_n)$ is a symplectic variety and has a unique 0-dimensional symplectic leaf $\{q\}$ in the sense of Theorem \ref{thm:Kaledin}.
\end{prop}

\begin{figure}[H]
\centering
\includegraphics[scale=0.3
]{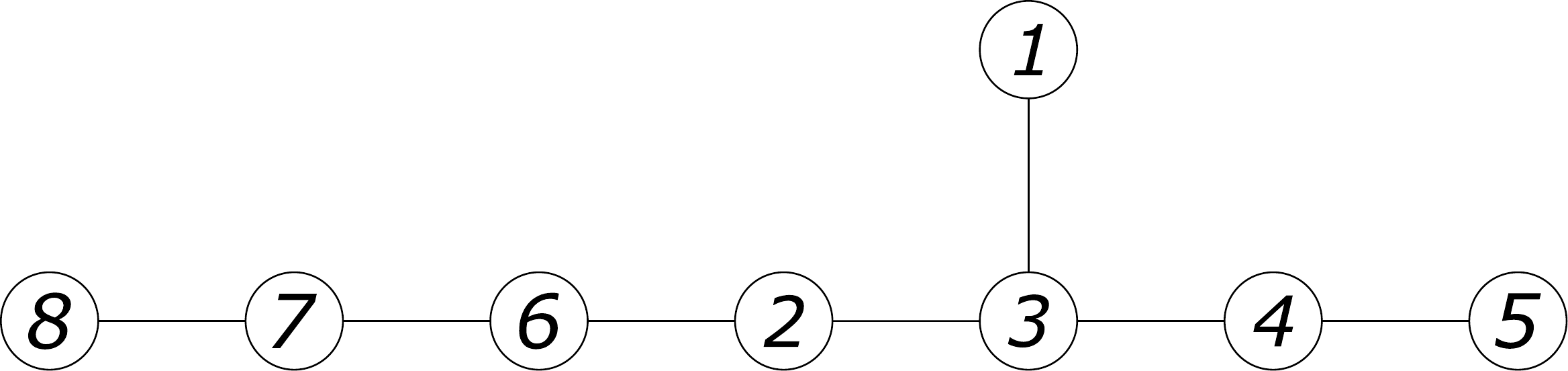}
\caption{The dual graph of the exceptional curves of type $\mathbf{E}_n$}
\label{fig6}
\end{figure}

Let $\tilde{\Gamma}_n$ be the minimal resolution of $\Gamma_n$ and let $\{C_1,\dots,C_n\}$ be the set of the exceptional curves ordered so that $C_i$ corresponds to the $i$-th vertex in the dual graph in Figure \ref{fig6}. Here we regard the graphs of $\mathbf{E}_6$ and $\mathbf{E}_7$ as subgraphs of $\mathbf{E}_8$.

We define the surfaces $P_{i,j}$ and $Q_i$ as before. One can directly check that successive blowing-ups along singular loci of $\mathrm{Hilb}^2(\Gamma_n)$ give a resolution and that this is symplectic using the same argument as Proposition \ref{prop:resolA}. This resolution is obtained from $\mathrm{Hilb}^2(\tilde{\Gamma}_n)$ by performing Mukai flops in the following order

{\footnotesize
$$\begin{aligned}
&(\text{the same as }\mathbf{D}_5)\to P_{6,6}\to P_{2,6}\to P_{3,6}\to P_{1,6}\to P_{4,6}\to P_{2,3}(2)\to P_{2,2}(2)\to Q_2\\
\to &P_{5,6}\to P_{2,4}(2)\to P_{3,3}(3)\to P_{1,4}(2)\to P_{1,2}(2)\to P_{1,3}(2)\to P_{1,1}(2)\to Q_1\,(\mathbf{E}_6 \text{ ends here})\\
\to &P_{7,7}\to P_{6,7}\to P_{2,7}\to P_{3,7}\to P_{1,7}\to P_{4,7}\to P_{3,6}(2)\to P_{2,6}(2)\to P_{6,6}(2)\to P_{5,7}\\
\to &P_{4,6}(2)\to P_{2,3}(3)\to P_{2,2}(3)\to P_{1,6}(2)\to P_{3,3}(4)\to P_{2,4}(3)\to P_{5,6}(2)\to Q_6\\
\to &Q_3(2)\to P_{4,4}(3)\to P_{3,4}(3)\to P_{2,5}(2)\to P_{1,5}(2)\to P_{3,5}(2)\to P_{4,5}(2)\to P_{5,5}(2)\\
\to &Q_5\,(\mathbf{E}_7 \text{ ends here})\to P_{8,8}\to P_{7,8}\to P_{6,8}\to P_{2,8}\to P_{3,8}\to P_{1,8}\to P_{4,8}\to P_{3,7}(2)\\
\to &P_{2,7}(2)\to P_{6,7}(2)\to P_{5,8}\to P_{4,7}(2)\to P_{3,6}(3)\to P_{2,6}(3)\to P_{1,7}(2)\to P_{2,3}(4)\\
\to &P_{4,6}(3)\to P_{5,7}(2)\to P_{7,7}(2)\to P_{6,6}(3)\to P_{2,2}(4)\to P_{3,3}(5)\to P_{2,4}(4)\to P_{1,6}(3)\\
\to &P_{3,4}(4)\to P_{4,4}(4)\to Q_3(3)\to Q_6(2)\to Q_7\to Q_2(2) \to P_{1,4}(3)\to P_{1,3}(3)\to P_{1,1}(3)\\
\to &P_{1,2}(3)\to P_{3,3}(6)\to P_{2,2}(5)\to P_{6,6}(4)\to P_{7,7}(3)\to Q_4(2)\to P_{4,6}(4)\to P_{2,3}(5)\\
\to &P_{1,7}(2)\to P_{2,6}(4)\to P_{3,6}(4)\to P_{4,7}(3)\to P_{5,8}(2)\to P_{6,7}(3)\to P_{2,7}(3)\to P_{3,7}(3)\\
\to &P_{4,8}(2)\to P_{1,8}(2)\to P_{3,8}(2)\to P_{2,8}(2)\to P_{6,8}(2)\to P_{7,8}(2)\to P_{8,8}(2)\to Q_8
\end{aligned}$$}
where $\mathbf{D}_5$ is also considered as the subgraph of $\mathbf{E}_8$. This order corresponds to the total order of positive roots for $\mathbf{E}_n$ defined in the same way as $\mathbf{D}_n$. The number of flops we need for $\mathbf{E}_6,\mathbf{E}_7$ and $\mathbf{E}_8$ are 36, 63, and 120 respectively.

\begin{figure}[p]
\centering
\includegraphics[scale=0.7]{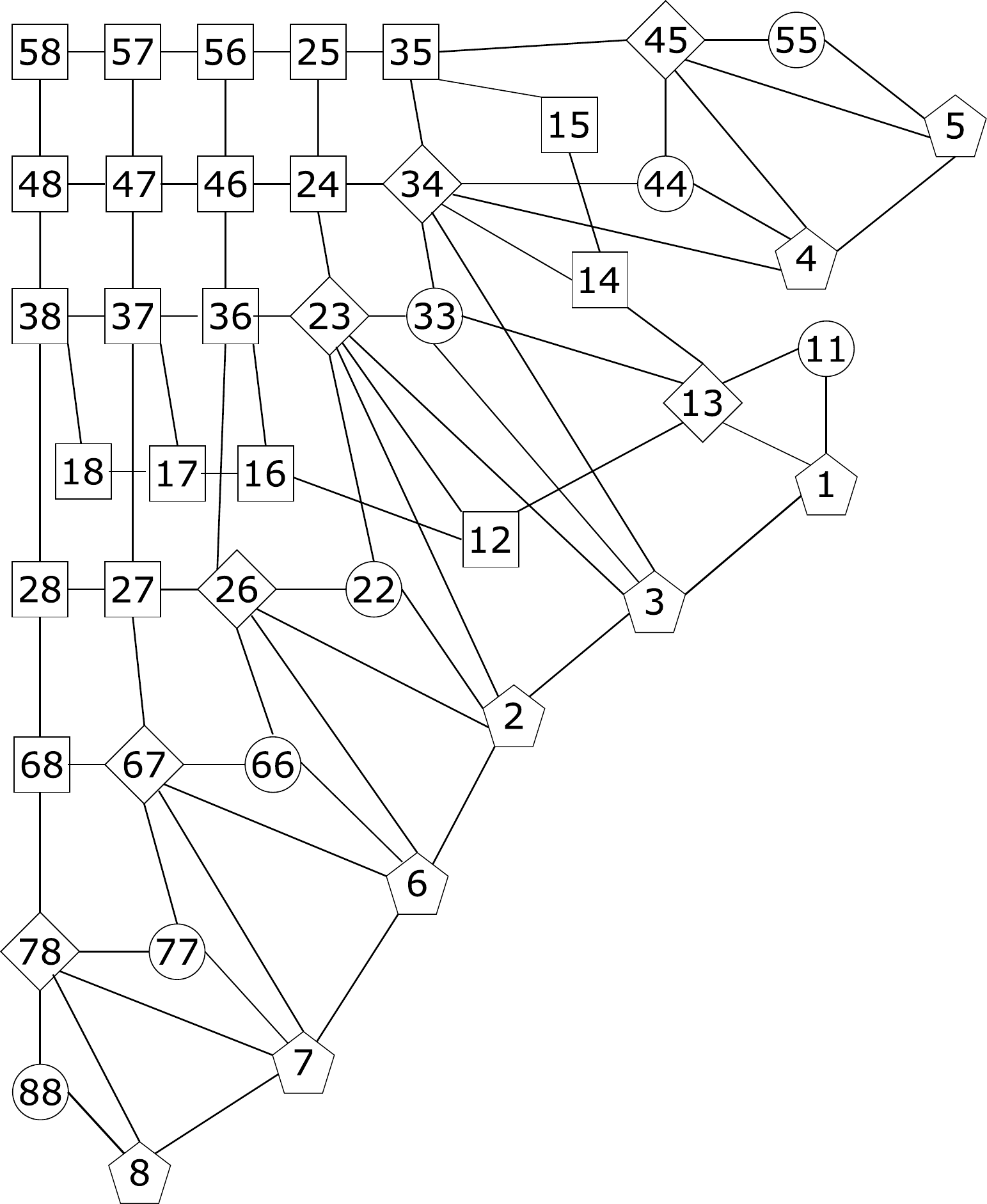}
\caption{The fiber of the Hilbert-Chow resolution $\pi_2\circ\pi_1$}
\label{fig7}
\end{figure}

\begin{figure}[p]
\centering
\includegraphics[scale=0.7]{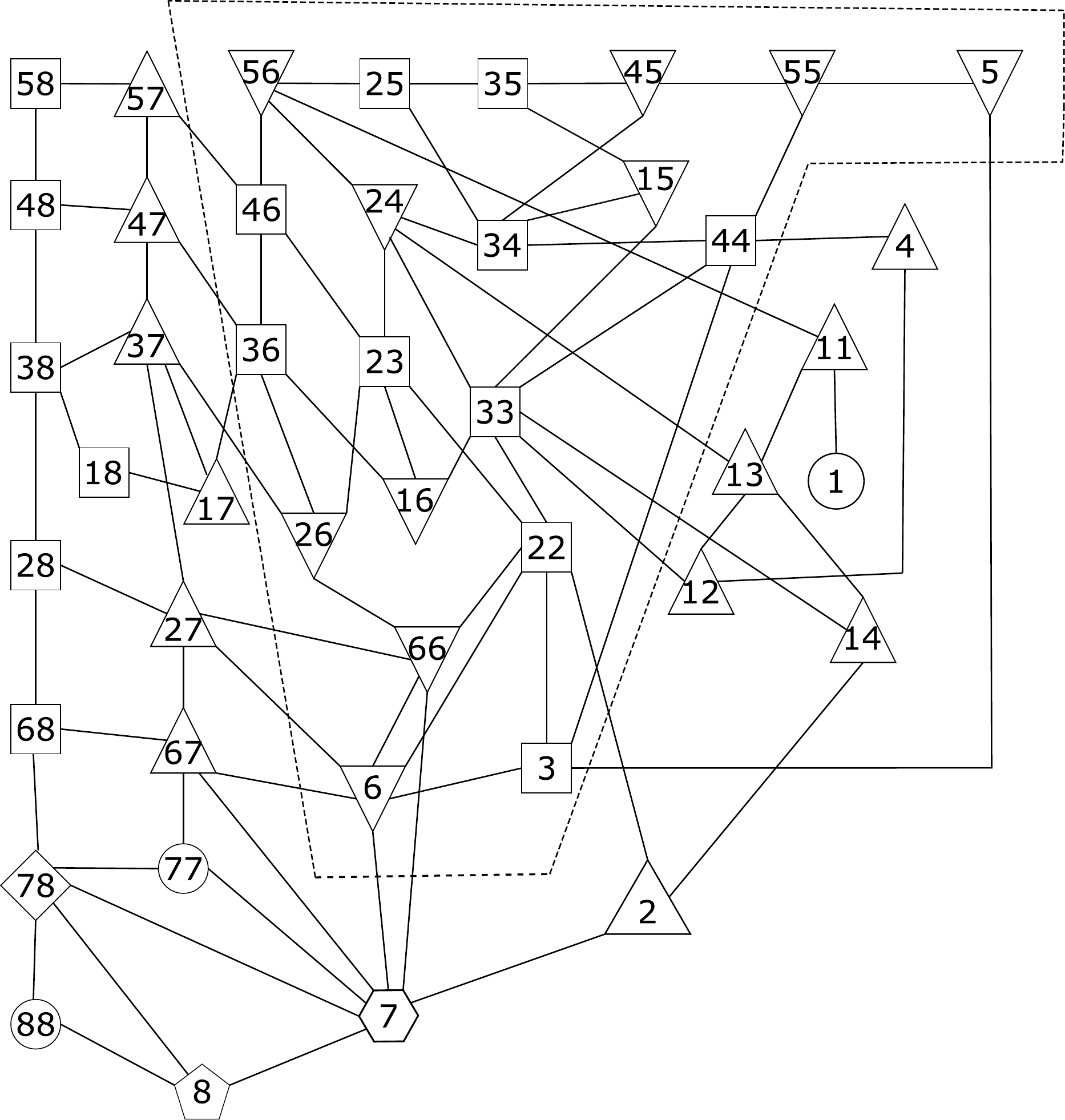}
\caption{The central fiber of $\psi$ for $\mathbf{E}_6$}
\label{fig8}
\end{figure}

\begin{figure}[p]
\centering
\includegraphics[scale=0.7]{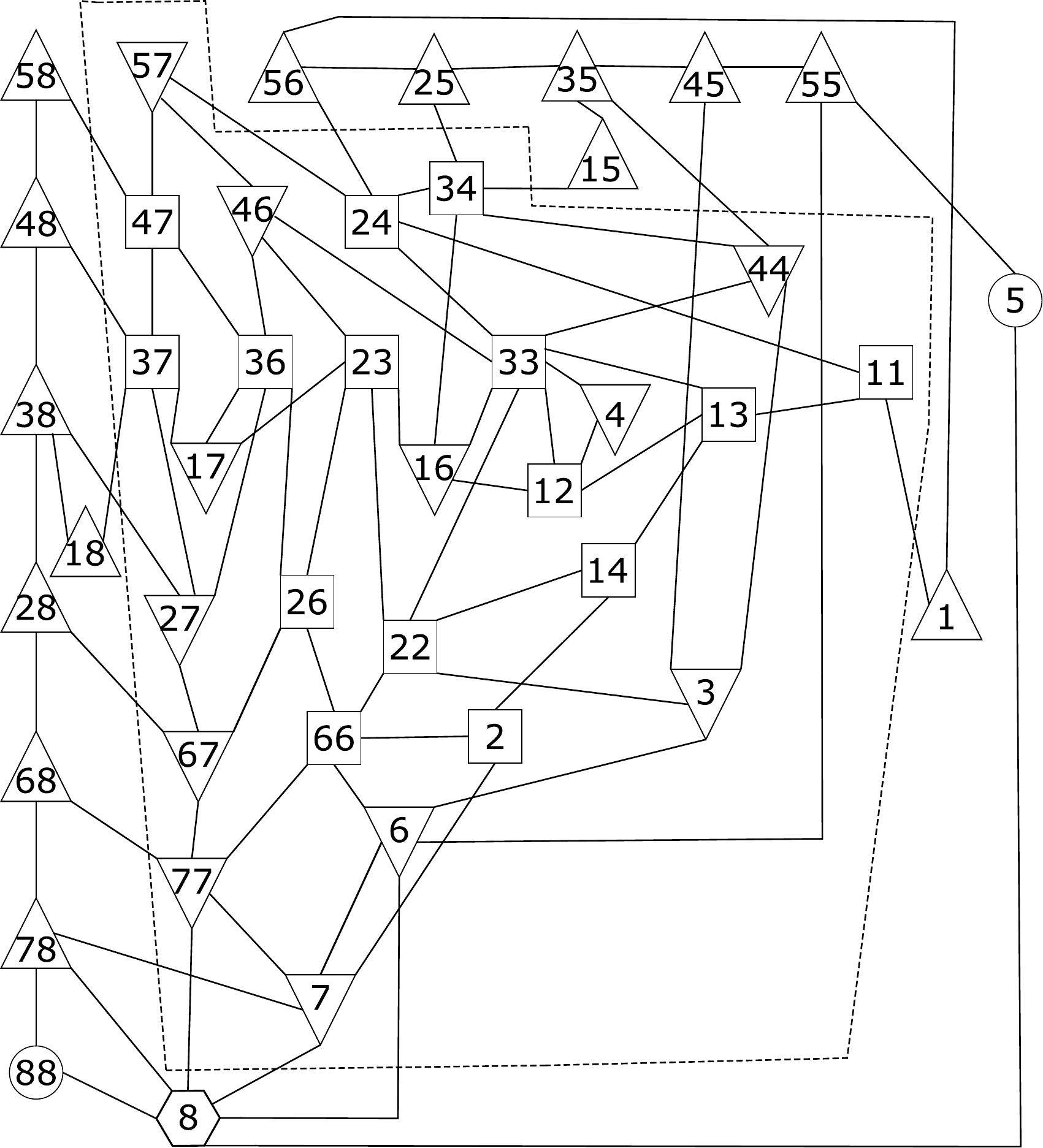}
\caption{The central fiber of $\psi$ for $\mathbf{E}_7$}
\label{fig9}
\end{figure}

\begin{figure}[p]
\centering
\includegraphics[scale=0.7]{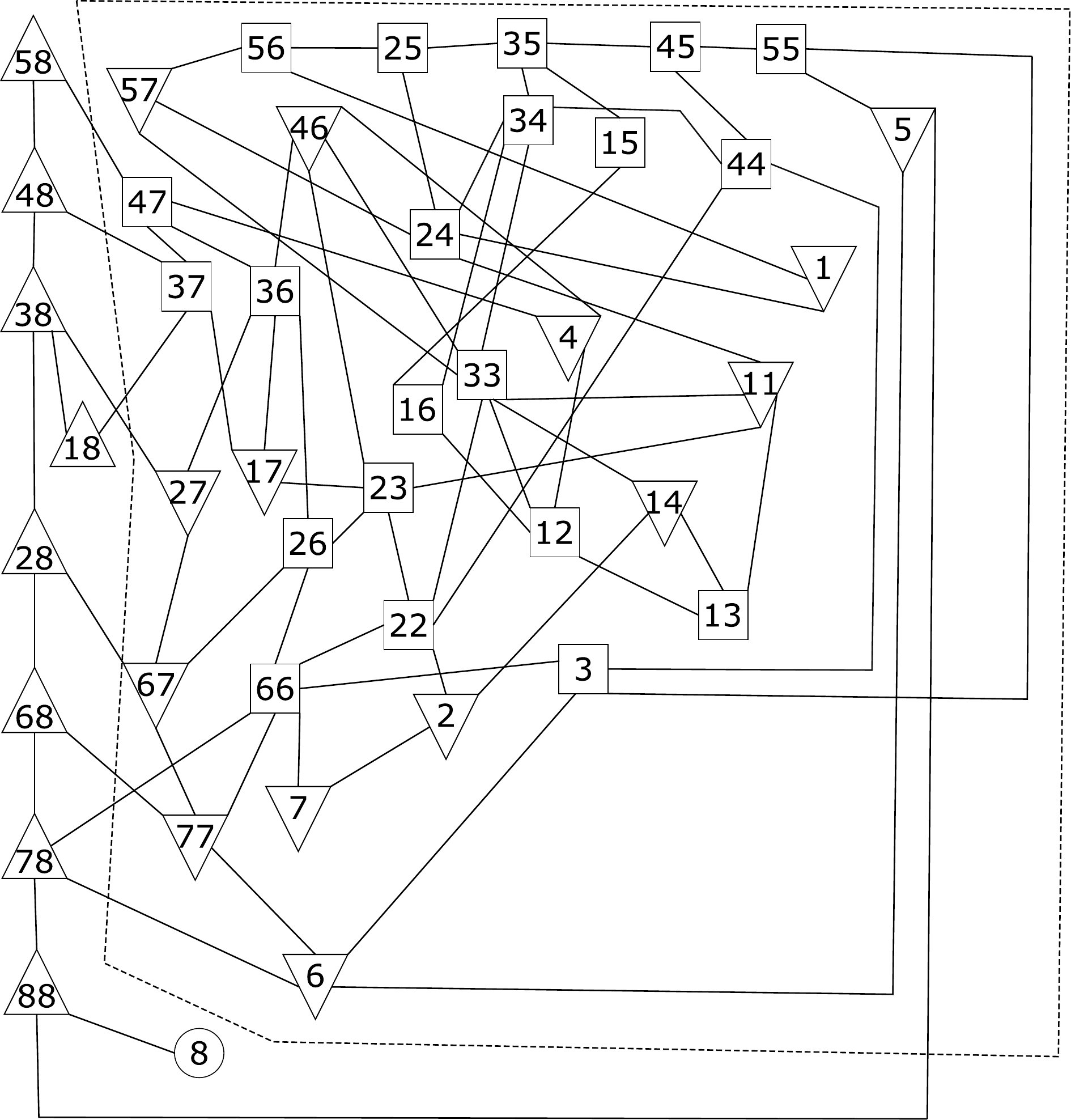}
\caption{The central fiber of and $\psi$  for $\mathbf{E}_8$}
\label{fig10}
\end{figure}

We give the initial diagram and the final diagram of the central fiber of the resolutions (Figure \ref{fig7}, \ref{fig8}, \ref{fig9}, and \ref{fig10}). We only treat the case $n=8$ since the diagrams for $\mathbf{E}_6$ and $\mathbf{E}_7$ are embedded in the $\mathbf{E}_8$-diagram in an obvious way. The presentation in the figures are the same as the case for $\mathbf{D}_n$.


The description of the central fibers shows that the symplectic resolutions are unique (see Remark \ref{rem2}). We give the intersections of the central fiber $F$ of $\psi$ and the exceptional divisors $E_k$ as well as types $\mathbf{A}_n$  and $\mathbf{D}_n$. Note that, unlike $\mathbf{A}_n\,(n\ge3)$ and $\mathbf{D}_n\,(n\ge5)$, the isomorphism classes of the components in $F$ can change from $\mathbf{E}_6$ to $\mathbf{E}_7$ or $\mathbf{E}_8$, though in fact the fibers $F$ for $\mathbf{E}_6$ and $\mathbf{E}_7$ have abstract embeddings into $F$ for $\mathbf{E}_8$ (see Section \ref{4}).

When $n=6$, $F_k=F\cap E_k$ is given as follows.
{\footnotesize
$$F_1=P_{2,5}\cup P_{3,5}\cup P_{3,6}\cup P_{4,6}\cup P_{5,6},\;
F_2=P_{1,5}\cup P_{2,2}\cup P_{3,3}\cup P_{3,5}\cup P_{6,6}$$
$$F_3=P_{2,3}\cup P_{2,4}\cup P_{2,6}\cup P_{3,4}\cup P_{4,5},\;
F_4=P_{1,6}\cup P_{3,3}\cup P_{3,6}\cup P_{4,4}\cup P_{5,5}$$
$$F_5=P_{2,2}\cup P_{2,3}\cup P_{4,6}\cup Q_3\cup Q_5,\;
F_6=P_{2,5}\cup P_{3,4}\cup P_{4,4}\cup Q_3\cup Q_6$$}

Each $F_k$ is a $\mathbb{P}^1$-bundle over a Dynkin tree of $\mathbb{P}^1$'s of type $\mathbf{A}_5$.

When $n=7$, we have:
{\footnotesize
$$F_1=P_{1,1}\cup P_{2,4}\cup P_{3,4}\cup P_{3,7}\cup P_{4,7}\cup P_{5,7},\;
F_2=P_{1,2}\cup P_{1,6}\cup P_{2,3}\cup P_{2,6}\cup P_{3,4}\cup P_{6,7}$$
$$F_3=P_{1,3}\cup P_{2,7}\cup P_{3,3}\cup P_{3,6}\cup P_{4,4}\cup P_{4,6},\;
F_4=P_{1,4}\cup P_{1,7}\cup P_{2,2}\cup P_{2,3}\cup P_{3,7}\cup Q_3$$
$$F_5=P_{2,6}\cup P_{3,6}\cup P_{4,7}\cup P_{6,6}\cup Q_2\cup Q_6,\;
F_6=P_{2,2}\cup P_{2,4}\cup P_{3,3}\cup P_{6,6}\cup P_{7,7}\cup Q_4$$
$$F_7=P_{1,1}\cup P_{1,2}\cup P_{1,3}\cup P_{1,4}\cup Q_2\cup Q_7.$$}
Each $F_k$ is a $\mathbb{P}^1$-bundle over a Dynkin tree of $\mathbb{P}^1$'s of type $\mathbf{D}_6$.

When $n=8$, we have:
{\footnotesize
$$ F_1=P_{1,2}\cup P_{1,3}\cup P_{1,5}\cup P_{1,6}\cup P_{3,7}\cup P_{4,7}\cup Q_4,\;
F_2=P_{1,1}\cup P_{1,3}\cup P_{2,3}\cup P_{2,4}\cup P_{2,5}\cup P_{2,6}\cup P_{6,7}$$
$$F_3=P_{1,4}\cup P_{2,7}\cup P_{3,3}\cup P_{3,4}\cup P_{3,5}\cup P_{3,6}\cup P_{4,6},\;
F_4=P_{1,7}\cup P_{2,2}\cup P_{2,3}\cup P_{3,7}\cup P_{4,4}\cup P_{4,5}\cup Q_2$$
$$F_5=P_{2,6}\cup P_{3,6}\cup P_{4,7}\cup P_{5,5}\cup P_{6,6}\cup Q_3\cup Q_7,\;
F_6=P_{1,2}\cup P_{2,2}\cup P_{3,3}\cup P_{5,6}\cup P_{5,7}\cup P_{6,6}\cup P_{7,7}$$
$$F_7=P_{1,6}\cup P_{2,4}\cup P_{3,4}\cup P_{4,4}\cup Q_1\cup Q_3\cup Q_6,\;
F_8=P_{1,5}\cup P_{2,5}\cup P_{3,5}\cup P_{4,5}\cup P_{5,5}\cup P_{5,6}\cup Q_5.$$}
Each $F_k$ is a $\mathbb{P}^1$-bundle over a Dynkin tree of $\mathbb{P}^1$'s of type $\mathbf{E}_7$.

Similarly to the cases of $\mathbf{A}_n$ and $\mathbf{D}_n$, these $\mathbb{P}^1$-bundle structures and the data of 1-dimensional intersections determine the isomorphism class of $F$.

\section{Formal neighborhoods of central fibers}\label{3}

In this section we show that the isomorphism class of the central fiber of the unique symplectic resolution of $\mathrm{Hilb}^2(\Gamma)$, which was described in the previous section, determines the local singularity (Theorem \ref{thm:main}). This is done by studying the formal neighborhoods of the central fiber.

Let $X$ be a 4-dimensional symplectic variety and $\pi:Y\to X$ a symplectic resolution. If $\mathrm{Sing}(X)$ is smooth, the germ of $X$ at any $p\in\mathrm{Sing}(X)$ is the same as that of $\mathbb{C}^2\times (\text{ADE-singularity})$ (see Theorem \ref{thm:Kaledin}). Therefore, we assume $X$ has a 0-dimensional symplectic leaf $\{p\}$.

It is known that each irreducible component $V$ of the central fiber $F=\pi^{-1}(p)$ is a Lagrangian subvariety of $Y$, i.e., $V$ is 2-dimensional and the restriction of the symplectic form on $Y$ to $V$ is zero (cf. \cite[Prop. 5.3.]{AW}). If $X$ is $\mathrm{Hilb}^2(\Gamma)$, then $V$ is isomorphic to $\Sigma_0=\mathbb{P}^1\times\mathbb{P}^1$ or $\Sigma_2$ as we saw in the previous section. One easily sees that, for any Lagrangian submanifold $V$ in $Y$, the symplectic form $T_Y\times T_Y\to \mathcal{O}_Y$ induces a non-degenerate pairing $T_V\times N_{V/Y}\to \mathcal{O}_V$ and thus the normal bundle of $V$ in $Y$ is isomorphic to the cotangent bundle.

The symplectic form on $Y$ can be regarded as an $\mathcal{O}_Y$-linear map $\Theta:\wedge^2 \Omega_Y^1\to \mathcal{O}_Y$, and it induces a non-degenerate Poisson structure $\{-,-\}:\mathcal{O}_Y\times\mathcal{O}_Y\to \mathcal{O}_Y$ defined by $\{f,g\}=\Theta(df\wedge dg)$ (see \cite[\S1]{K}). For any ideal $\mathcal{I}\subset\mathcal{O}_Y$ and $n\in\mathbb{N}$, the Poisson structure of $Y$ naturally induces a $\mathbb{C}$-skew-bilinear map $\mathcal{O}_Y/\mathcal{I}^n\times\mathcal{O}_Y/\mathcal{I}^n\to \mathcal{O}_Y/\mathcal{I}^{n-1}$ that has the same properties as the Poisson bracket i.e., the Leibniz rule and the Jacobi identity. (By abuse of language, we also call such a bilinear map {\em a Poisson structure}.) In particular, the formal completion of $\mathcal{O}_Y$ along $\mathcal{I}$ has a natural Poisson structure.

When $\mathcal{I}$ defines a Lagrangian submanifold $V\subset Y$, the symplectic form gives an identification $T_V\cong\mathcal{I}/\mathcal{I}^2$. From the definition of the Poisson structure, one sees that the restriction of the map $\mathcal{O}_Y/\mathcal{I}^{n+1}\times\mathcal{O}_Y/\mathcal{I}^{n+1}\to \mathcal{O}_Y/\mathcal{I}^n$ to $\mathcal{I}^n/\mathcal{I}^{n+1}\times\mathcal{O}_Y/\mathcal{I}^{n+1}$ can be identified with the map
$$\mathrm{Sym}^n T_V\times \mathcal{O}_Y/\mathcal{I}^{n+1}\to\mathrm{Sym}^{n-1} T_V\subset \mathcal{O}_Y/\mathcal{I}^n,\,(\partial,x)\mapsto \partial(\bar{x})$$
where $\bar{x}$ denotes the image of $x$ in $\mathcal{O}_Y/\mathcal{I}$ and $\mathrm{Sym}^n T_V$ acts on $\mathcal{O}_Y/\mathcal{I}=\mathcal{O}_V$ as the natural derivation.

For any smooth variety $U$, the total space of the cotangent bundle $T^*_U$ has a natural symplectic structure and thus also has a Poisson structure $\{-,-\}$. Note that the restriction of the bracket $\{-,-\}$ to $T_U\times \mathcal{O}_U$ (resp. $T_U\times T_U$) is the same as the derivation (resp. Lie bracket). Also, the restriction to $\mathcal{O}_U\times \mathcal{O}_U$ is zero. The following proposition shows that every symplectic manifold is locally a cotangent bundle in the formal neighborhood of a Lagrangian submanifold.

\begin{prop}\label{prop:local}
Let $Y$ be a nonsingular symplectic variety and $V\subset Y$ a nonsingular Lagrangian subvariety with the defining ideal $\mathcal{I}\subset \mathcal{O}_Y$. Then, for any closed point $p\in V$, there is an affine open neighborhood $U\subset V$ of $p$ such that the completion $\varprojlim\,(\mathcal{O}_{Y}/\mathcal{I}^n)|_U$ is isomorphic to $\prod_{n=0}^\infty \mathrm{Sym}^n T_U$ as a sheaf of Poisson algebras on $U$.
\end{prop}

{\em Proof.} Since $V$ is smooth, we can choose an affine open neighborhood $U\subset V$ of $p$ such that $\Omega^1_U$ is a free $\mathcal{O}_U$-module. We take $x_1,\dots,x_\ell\in H^0(\mathcal{O}_U)$ so that $dx_1,\dots,dx_\ell$ form a free basis of $H^0(\Omega^1_U)$. We will inductively construct compatible isomorphisms $\alpha_n:\prod_{k=0}^{n-1} \mathrm{Sym}^k T_U\to (\mathcal{O}_Y/\mathcal{I}^n)|_U,\,n=1,2,\dots$ such that the restriction of $\alpha_n$ to $\mathrm{Sym}^{n-1} T_U$ is the isomorphism $\mathrm{Sym}^{n-1} T_U\to (\mathcal{I}^{n-1}/\mathcal{I}^n)_U$ that is induced by the symplectic form on $Y$ and that $\alpha_n$'s preserve the Poisson structures i.e., $\{\alpha_{n+1}(x),\alpha_{n+1}(y)\}=\alpha_n(\{x,y\})$ for any sections $x,y\in\prod_{k=0}^n \mathrm{Sym}^k T_U$.

The case $n=1$ is trivial. So we assume that we have constructed the isomorphisms up to $\alpha_n$. Since $U$ is smooth and affine, the inclusion $\mathcal{O}_U\to(\mathcal{O}_{Y}/\mathcal{I}^n)|_U$ lifts to $(\mathcal{O}_{Y}/\mathcal{I}^{n+1})|_U$ by the infinitesimal lifting property. Also, the inclusion $(\mathcal{I}/\mathcal{I}^2)_U\to (\mathcal{O}_{Y}/\mathcal{I}^n)|_U$ lifts to $(\mathcal{O}_{Y}/\mathcal{I}^{n+1})|_U$ for each $k=1,\dots,n-1$ since $(\mathcal{I}/\mathcal{I}^2)_U$ is a free $\mathcal{O}_U$-module. Using the identification $T_U\cong(\mathcal{I}/\mathcal{I}^2)_U$, these liftings define an algebra isomorpism $\alpha'_{n+1}:\prod_{k=0}^n \mathrm{Sym}^k T_U\to (\mathcal{O}_{Y}/\mathcal{I}^{n+1})|_U$ but it possibly does not preserve the Poisson structures.

Let $\{-,-\}'$ be the Poisson structure on $\prod_{k=0}^n \mathrm{Sym}^k T_U$ which is the pullback of the Poisson structure on $(\mathcal{O}_{Y}/\mathcal{I}^{n+1})|_U$ under $\alpha'_{n+1}$. Note that $\{-,-\}'$ is non-degenerate i.e., for any $\mathbb{C}$-derivation $\partial$ on $\prod_{k=0}^n \mathrm{Sym}^k T_U$, there is $f\in \prod_{k=0}^n \mathrm{Sym}^k T_U$ such that $\pi(\partial(-))=\{f,-\}'$ where $\pi:\prod_{k=0}^n \mathrm{Sym}^k T_U\to\prod_{k=0}^{n-1} \mathrm{Sym}^k T_U$ is the natural projection.

Let $\partial_1,\dots,\partial_\ell\in T_U$ be the dual basis of $dx_1,\dots,dx_\ell$. Let $f_i\in \prod_{k=0}^n \mathrm{Sym}^k T_U$ be an element that corresponds to the derivation $\{x_i,-\}$ via the non-degenerate bracket $\{-,-\}'$. Since the two Poisson structures are compatible with $\pi$ by assumption, we may assume that the difference $v(x_i):=f_i-x_i$ is in $\mathrm{Sym}^n T_U$ by replacing $f_i$ by $f_i+c_i$ with a unique constant $c_i\in\mathbb{C}$. Then the assignment $x_i\mapsto v(x_i)$ defines a $\mathbb{C}$-derivation $v:\mathcal{O}_U\to \mathrm{Sym}^n T_U$. Note that $v(x_i)(x_j)=\{x_i,x_j\}-\{x_i,x_j\}'$. If one replaces $\alpha'_n$ by $\alpha'_n\circ(\mathrm{Id}+\frac{1}{2}v)$, which is still an algebra isomorphism, one sees that $\{\mathcal{O}_U,\mathcal{O}_U\}'=0$.

Next we take elements $g_1,\dots,g_\ell$ of $\prod_{k=0}^n \mathrm{Sym}^k T_U$ that correspond to the derivations $\{\partial_i,-\}$ using the non-degeneracy of $\{-,-\}'$. Then, similarly as above, the differences $w(\partial_i):=g_i-\partial_i$ define an $\mathcal{O}_U$-linear map $w:T_U\to \mathrm{Sym}^n T_U$, and replacing $\alpha'_n$ by $\alpha_n:=\alpha'_n\circ(\mathrm{Id}+w)$ gives a bracket $\{-,-\}'$ that satisfies $\{\partial_i,x_j\}'=\delta_{i,j}$ where $\delta_{i,j}$ is the Kronecker delta. Note that $w$ is a unique one satisfying this property, whereas the above $v$ is not. Using the non-degeneracy again, we see that $\{\partial_i,\partial_j\}'=0$ for any $i,j$. Therefore, we have $\{-,-\}'=\{-,-\}$, and this completes the proof.
\qed

\vspace{3mm}

Symplectic manifolds may not be cotangent bundles globally on the Lagrangian submanifold. The following lemma gives a cohomological criterion for this. 

\begin{lem}\label{lem:normal}
Let $V$ be a nonsingular variety. Suppose that we have the following exact sequence of sheaves on $V$ for each $l=0,1,\dots,n$:\\
$$0\to\mathrm{Sym}^l T_V\to \mathcal{A}_l\to \mathcal{A}_{l-1}\to0$$
satisfying\\
$\bullet$ $\mathcal{A}_{-1}=0$, and each $\mathcal{A}_l\,(l\ge0)$ is an infinitesimal deformation of $\mathcal{O}_V$ i.e. a sheaf of finitely generated $\mathbb{C}$-algebras on $V$ such that $\mathcal{A}_l/\mathrm{nil}(\mathcal{A}_l)\cong\mathcal{O}_V$ where $\mathrm{nil}(\mathcal{A}_l)$ is the nilradical of $\mathcal{A}_l$,\\
$\bullet$ the map $\mathcal{A}_l\to \mathcal{A}_{l-1}$ is an algebra surjection such that $\mathrm{Sym}^l T_V$ is square-zero as an ideal of $\mathcal{A}_l$,\\
$\bullet$ each $\mathcal{A}_l\,(l\ge0)$ admits a non-degenerate $\mathbb{C}$-bilinear Poisson structure $\{-,-\}:\mathcal{A}_l\times\mathcal{A}_l\to \mathcal{A}_{l-1}$ which is compatible with $\mathcal{A}_l\to \mathcal{A}_{l-1}$, and\\
$\bullet$ there is an open affine covering $V=\bigcup_i U_i$ such that we have compatible algebra isomorphisms $(\alpha_l)_i:\prod_{k=0}^l \mathrm{Sym}^k T_{U_i}\to \mathcal{A}_l|_{U_i}$ for $l=1,\dots,n-1$ which are the identities on $\mathrm{Sym}^l T_{U_i}$ and preserve the Poisson structures, i.e. $
\{(\alpha_l)_i(x),(\alpha_l)_i(y)\}=(\alpha_{l-1})_i(\{x,y\})\in \mathcal{A}_{l-1}|_{U_i}$ for any $x,y\in \prod_{k=0}^l T_{U_i}$.


If we are given two such collections $\mathcal{A}=\{(\mathcal{A}_l,(\alpha_l)_i)\}_{l=0,1,\dots,n},\,\mathcal{A}'=\{(\mathcal{A}'_l,(\alpha'_l)_i)\}_{l=0,1,\dots,n}$, and an algebra isomorphism $\alpha_{n-1}:\mathcal{A}_{n-1}\to\mathcal{A}'_{n-1}$ such that $\alpha_{n-1}|_{U_i}=(\alpha'_{n-1})_i\circ(\alpha_{n-1})_i^{-1}$. Then the pair $(\mathcal{A},\mathcal{A}')$ naturally defines a class $\xi$ in $H^1(V,T_V\otimes\mathrm{Sym}^n T_V)$ such that $\xi=0$ if and only if there is a Poisson isomorphism between the two exact sequences for $\mathcal{A}_n$ and $\mathcal{A}'_n$ (i.e. a Poisson algebra isomorphism $\alpha_n:\mathcal{A}_n\to\mathcal{A}'_n$ which is compatible with the maps in the sequences for $l=n$, the identity map of $\mathrm{Sym}^n T_V$, and $\alpha_{n-1}$).

In particular, if $V$ is a nonsingular Lagrangian subvariety of  a nonsingular symplectic variety $Y$ and $H^1(V,T_V\otimes\mathrm{Sym}^n T_V)$ vanish for all $n\ge1$, then $\hat{Y}_V$ is isomorphic as a formal Poisson scheme to the completion of the total space of the cotangent bundle $T^*_V$ along its zero section.
\end{lem}

{\em Proof.} The non-degeneracy of the Poisson structure implies that there is a lift $(\alpha_n)_i:\prod_{k=0}^n \mathrm{Sym}^k T_{U_i}\to \mathcal{A}_n|_{U_i}$ of $(\alpha_{n-1})_i$ which is the identity on  $\mathrm{Sym}^n T_{U_i}$ and satisfies $\{(\alpha_n)_i(\partial),(\alpha_n)_i(x)\}=\alpha_{n-1}(\partial(x))$ for any $\partial\in T_{U_i}$ and $x\in \mathcal{O}_{U_i}$. Indeed, we already know that $\prod_{k=0}^n \mathrm{Sym}^k T_{U_i}$ and $\mathcal{A}_n|_{U_i}$ are isomorphic  
as algebras by Proposition \ref{prop:local}, and we can take such $(\alpha_n)_i$ by replacing an isomorphism using the $\mathcal{O}_{U_i}$-linear map $w$ defined in the final paragraph of the proof of Proposition \ref{prop:local}. Note that we do not have to choose a suitable derivation $v:\mathcal{O}_{U_i}\to \mathrm{Sym}^n T_{U_i}$ as in the proof of the proposition and that $(\alpha_n)_i$ is not required to preserve the whole Poisson structures at this stage.

Set $U_{i,j}=U_i\cap U_j$. Then we obtain an automorphism $(\alpha_n)_{i,j}:=(\alpha_n)_i^{-1}\circ (\alpha_n)_j$ of $\prod_{k=0}^n \mathrm{Sym}^k T_{U_{i,j}}$ for each $i,j$. Similarly, we consider the corresponding automorphisms $(\alpha_n)'_{i,j}$ for $\mathcal{A}'_n$. Then $d_{i,j}:=(\alpha_n)_{i,j}-(\alpha_n)'_{i,j}$ is a $\mathbb{C}$-derivation $\prod_{k=0}^n \mathrm{Sym}^k T_{U_{i,j}}\to \mathrm{Sym}^n T_{U_{i,j}}$. Such $d_{i,j}$ is regarded as an element of the vector space
$$\mathrm{Hom}_{\mathcal{O}_{U_{i,j}}}(\Omega_{U_{i,j}}^1,\mathrm{Sym}^n T_{U_{i,j}})\oplus\mathrm{Hom}_{\mathcal{O}_{U_{i,j}}}(T_{U_{i,j}},\mathrm{Sym}^n T_{U_{i,j}}).$$
Since both $(\alpha_n)_{i,j}$ and $(\alpha_n)'_{i,j}$ preserve the Poisson brackets on $T_{U_{i,j}}\times \mathcal{O}_{U_{i,j}}$, we see that the component of $d_{i,j}$ in the second direct summand should be zero by the fact that the restriction of the bracket to $\mathrm{Sym}^n T_{U_{i,j}}\times \mathcal{O}_{U_{i,j}}$ is the same as the natural derivation. (Note that this is equivalent to the uniqueness of the above $w$.)

One can check that $\{d_{i,j}\}_{i,j}$ is a \v{C}ech 1-cocycle of $T_V\otimes\mathrm{Sym}^n T_V$ and thus defines an element of $H^1(V,T_V\otimes\mathrm{Sym}^n T_V)$. One can also check that this is zero if and only if there is a possibly non-Poisson isomorphism between the two sequences.

When the cohomology class is zero, we can always take an isomorphism as a Poisson one. This is shown as follows. Recall that this is locally done by replacing the isomorphism using the derivation $v:\mathcal{O}_{U_i}\to \mathrm{Sym}^n T_{U_i}$ for each $n$ in the proof of Proposition \ref{prop:local}. One can check that in fact $v$ extends to a well-defined global derivation $\mathcal{O}_V\to \mathrm{Sym}^n T_V$ since $v$ can be defined just using the non-degeneracy of the Poisson bracket without using local coordinates. This proves the first claim.

For the second claim, let $\mathcal{I}\subset\mathcal{O}_Y$ be the defining ideal of $V$. Then, by Proposition \ref{prop:local}, the collection $\{\mathcal{O}_Y/\mathcal{I}^{n+1}\}_{n=-1,0,\dots}$ satisfies the assumptions on $\mathcal{A}_n$ where $\mathrm{Sym}^n T_V$ is identified with $\mathcal{I}^n/\mathcal{I}^{n+1}\subset\mathcal{O}_Y/\mathcal{I}^{n+1}$ using the symplectic form on $Y$. Assume that $H^1(V,T_V\otimes\mathrm{Sym}^n T_V)=0$. Then applying the first claim inductively to $\mathcal{A}_n=\mathcal{O}_Y/\mathcal{I}^{n+1}$ and $\mathcal{A}'_n=\prod_{k=0}^n \mathrm{Sym}^k T_V$ shows that there is a Poisson algebra isomorphism between $\varprojlim\,\mathcal{O}_Y/\mathcal{I}^{n+1}$ and $\prod_{n=0}^\infty \mathrm{Sym}^n T_V$.
\qed

\vspace{3mm}

In the case when $V$ is isomorphic to $\Sigma_0$ in a symplectic 4-fold $Y$, one can check that all the $H^1$'s vanish and thus the formal neighborhood $\hat{Y}_V$ of $V$ in $Y$ is isomorphic to the cotangent bundle of $V$. On the other hand, if $V\cong\Sigma_2$, the $H^1$'s do not vanish. The formal neighborhood of $V$ in the symplectic resolution of $\mathrm{Hilb}^2(\Gamma)$ is actually not isomorphic to the cotangent bundle as we will see in Remark \ref{rem3}. However, we can show that the isomorphism class of the formal neighborhood of $V=\Sigma_2$ is uniquely determined by that of the negative section of $\Sigma_2$.

\begin{lem}\label{lem:S2}
Let $Y$ be a nonsingular symplectic 4-fold and $V\subset Y$ a closed subvariety which is isomorphic to $\Sigma_2$. Let $s\subset V$ be the unique negative section. Assume that the first order thickening of $s$ in $Y$ is trivial i.e. isomorphic to the first order thickening in the normal bundle. Then the Poisson isomorphism class of the formal completion $\hat{Y}_V$ of $Y$ along $V$ is uniquely determined (independent of the embedding $V\subset Y$) by the isomorphism class of the formal completion $\hat{Y}_s$ of $Y$ along $s$.
\end{lem}

{\em Proof.} First note that $H^0(\Sigma_2,\Omega_{\Sigma_2}^2)=0$ and thus $V\subset Y$ is Lagrangian. Assume that we are given two embeddings $V\subset Y$ and $V\subset Y'$ and that the completions $\hat{Y}_s$ and $\hat{Y}'_s$ are isomorphic. We will show the vanishings of the cohomology classes $\xi_n\in H^1(V,T_V\otimes\mathrm{Sym}^n T_V)$ defined by the pair $Y$ and $Y'$ (see Lemma \ref{lem:normal}).


Let $T_n$ be the completion of $T_V\otimes\mathrm{Sym}^n T_V$ along $s$. Then $H^1(s,T_1)$ contains $H^1(s, T_s\otimes T_V|_s)$ as a direct summand. Let $\eta$ be the corresponding component of the restriction $\xi_1|_{\hat{V}_s}$. By the  triviality assumption, the first order thickenings of $s$ in $Y$ and $Y'$ are isomorphic as infinitesimal deformations of $s_2$ where $s_k$ is the $(k-1)$-th order thickening of $s$ in $V$. Therefore, we have $\eta=0$.

One sees that the sum of the restrictions of $\xi_n\,(n\ge 2)$ and $\xi_1-\eta$ to $\hat{V}_s$ sits inside the $\mathbb{C}$-vector space
$$\prod_{k\ge1}(H^1(s, T_s\otimes \mathrm{Sym}^{k+1} N^*)\oplus H^1(s, N\otimes\mathrm{Sym}^k N^*))$$
where $N\cong\mathcal{O}_s(2)\oplus\mathcal{O}_s(-2)^{\oplus2}$ is the normal bundle of $s$ in $Y$. Using the argument of the proof of Lemma \ref{lem:normal}, this vector space is considered as the obstruction space for the two infinitesimal deformations $\mathcal{O}_{\hat{Y}_s}$ and $\mathcal{O}_{\hat{Y}'_s}$ of $\mathcal{O}_s$ to be isomorphic (cf. \cite[\S4, Satz 7.]{Gr}). Thus the cohomology class in this space is zero using the assumption that $\hat{Y}_s$ and $\hat{Y}'_s$ are isomorphic.


Therefore, we have to show that the restriction maps $H^1(V,T_V\otimes\mathrm{Sym}^n T_V)\to H^1(s,T_n)$ are injective for all $n$. To prove the claim, it suffices to show that the maps $H^1(V,\mathcal{F})\to H^1(s_k,\mathcal{F}|_{s_k})$ are injective for $\mathcal{F}=T_V\otimes\mathrm{Sym}^n\,T_V$ and $k\gg0$. We first show that this is true when $\mathcal{F}$ is an effective line bundle on $V$. Let $f\subset V$ be a fiber of $V\to\mathbb{P}^1$.  Then every effective line bundle is of the form $L=\mathcal{O}_V(as+bf)$ with $a,b\ge0$. Using the short exact sequence
$$0\to L(-ks)\to L\to L|_{s_k}\to0$$ 
we see that $H^1(V,L)\to H^1(s_k,L|_{s_k})$ is surjective since $H^2(V,L(-ks))\cong H^0(V,L(ks)\otimes \omega_V)^\vee=H^0(V,\mathcal{O}_V((k-a-2)s-(b+4)f))^\vee=0$. Applying the Leray spectral sequence to the retraction $s_k\to s$, we see that, for sufficiently large $k$,
$$\begin{aligned}h^1(s_k,L|_{s_k})&=h^1(s,(\bigoplus_{r=0}^{k+1}\mathcal{O}_s(2r))\otimes L|_s)\\
&=h^1(s,\bigoplus_{r=0}^{k+1}\mathcal{O}_s(2r-2a+b))=\sum_{r=0}^{a-[b/2]-1}(-2r+2a-b-1).
\end{aligned}$$
It is not hard to check that $h^1(V,L)$ is the same as this number. Therefore, $H^1(V,L)\to H^1(s_k,L|_{s_k})$ is an isomorphism for $k\gg0$.

Using the generalized Euler sequence \cite[Thm. 8.1.6]{CLS}
$$0\to\mathcal{O}_V^2\to \mathcal{O}_V(f)^2\oplus\mathcal{O}_V(s)\oplus\mathcal{O}_V(s+2f)\to T_V\to0,$$
we can prove the case when $\mathcal{F}=T_V$. For the symmetric powers of $T_V$, we may use the filtration (see \cite[Ch. 2. Ex. 5.16(c)]{Ha}) associated to the generalized Euler sequence. Then we can finally conclude that $H^1(V,\mathcal{F})\to H^1(s_k,\mathcal{F}|_{s_k})$ is an isomorphism for $\mathcal{F}=T_V\otimes\mathrm{Sym}^n\,T_V$ with $k\gg0$.
\qed

\vspace{3mm}

To determine the formal neighborhood of the fiber $F$ of $\mathrm{Hilb}^2(\Gamma)$, we will need to glue the formal neighborhoods of irreducible components. This is justified by the following lemma.

\begin{lem}\label{lem:glue}
Let $Y_1$ and $Y_2$ be algebraic schemes (over $\mathbb{C}$) and $V_i,W_i\subset Y_i,\,i=1,2$ be closed subschemes. Assume that there are isomorphisms of the formal completions $\phi_V:{(\hat{Y_1})}_{V_1}\to {(\hat{Y_2})}_{V_2}$ and $\phi_W:{(\hat{Y_1})}_{W_1}\to {(\hat{Y_2})}_{W_2}$ such that their restrictions to ${(\hat{Y_1})}_{V_1\cap W_1}$ induce the same isomorphism $\phi:{(\hat{Y_1})}_{V_1\cap W_1}\to {(\hat{Y_2})}_{V_2\cap W_2}$. Then there is an isomorphism of ${(\hat{Y_1})}_{V_1\cup W_1}$ and ${(\hat{Y_2})}_{V_2\cup W_2}$ which is restricted to $\phi_V$ and $\phi_W$ where the unions of the subschemes are defined by the intersections of the corresponding ideal sheaves.
\end{lem}

{\em Proof.} Let $\mathcal{I}_i,\mathcal{J}_i\subset\mathcal{O}_{Y_i},\,i=1,2$ be the ideal sheaves of $V_i$ and $W_i$. By assumption we have isomorphisms of the completions
$$\psi_V:\varprojlim\,\mathcal{O}_{Y_1}/\mathcal{I}_1^n\cong\varprojlim\,\mathcal{O}_{Y_2}/\mathcal{I}_2^n$$
and
$$\psi_W:\varprojlim\,\mathcal{O}_{Y_1}/\mathcal{J}_1^n\cong\varprojlim\,\mathcal{O}_{Y_2}/\mathcal{J}_2^n$$
which are restricted to the same $\varprojlim\,\mathcal{O}_{Y_1}/(\mathcal{I}_1+\mathcal{J}_1)^n
\cong\varprojlim\,\mathcal{O}_{Y_2}/(\mathcal{I}_2+\mathcal{J}_2)^n$. We must show that there is an isomorphism $\psi:\varprojlim\,\mathcal{O}_{Y_1}/(\mathcal{I}_1\cap\mathcal{J}_1)^n
\cong\varprojlim\,\mathcal{O}_{Y_2}/(\mathcal{I}_2\cap\mathcal{J}_2)^n$ compatible with $\psi_V$ and $\psi_W$. We will construct a homomorphism
$$\varprojlim\,\mathcal{O}_{Y_1}/(\mathcal{I}_1\cap\mathcal{J}_1)^n
\to\mathcal{O}_{Y_2}/(\mathcal{I}_2\cap\mathcal{J}_2)^m$$
for each $m\in\mathbb{N}$.

One can see that $(\mathcal{I}_2\cap\mathcal{J}_2)^m
\supset\mathcal{I}_2^n\cap\mathcal{J}_2^n$ for $n\gg m$ by considering the primary decompositions of the both sides. Thus, it suffices to construct a map to $\mathcal{O}_{Y_2}/(\mathcal{I}_2^n\cap\mathcal{J}_2^n)$. Let $x$ be any local section of $\varprojlim\,\mathcal{O}_{Y_1}/(\mathcal{I}_1\cap\mathcal{J}_1)^n$. Then we obtain local sections $x_1$ and $x_2$ of $\mathcal{O}_{Y_2}/\mathcal{I}_2^n$ and $\mathcal{O}_{Y_2}/\mathcal{J}_2^n$ respectively by using $\psi_V$ and $\psi_W$. The local sections $x_1$ and $x_2$ are mapped to the same local section of $\mathcal{O}_{Y_2}/(\mathcal{I}_2^n
+\mathcal{J}_2^n)$ by the compatibility of $\psi_V$ and $\psi_W$ and by the fact that $\mathcal{I}_2^n+\mathcal{J}_2^n
\supset(\mathcal{I}_2+\mathcal{J}_2)^k$ for $k\gg n$. Therefore, there exits a unique local section of $\mathcal{O}_{Y_2}/(\mathcal{I}_2^n\cap\mathcal{J}_2^n)$ mapped to $x_1$ and $x_2$. This construction is canonical in the sense that it gives a desired well-defined map $\psi:\varprojlim\,\mathcal{O}_{Y_1}/(\mathcal{I}_1\cap\mathcal{J}_1)^n
\to\varprojlim\,\mathcal{O}_{Y_2}/(\mathcal{I}_2\cap\mathcal{J}_2)^n$. The inverse of $\psi$ is constructed in the same way.
\qed

\vspace{3mm}

Let $Y$ be a symplectic resolution of (a germ of) a symplectic variety such that its central fiber $F$ is isomorphic to that of $\mathrm{Hilb}^2(\Gamma)$. We will show that the isomorphism class of the formal neighborhood of $F\subset Y$ is uniquely determined from $F$ itself. 

We first consider the simplest case where $\Gamma$ is of type $\mathbf{A}_3$. In this case $F$ is a union of two components $V\cong\Sigma_0$ and $W\cong\Sigma_2$. The intersection $s=V\cap W$ is a diagonal of $V$ and the negative section of $W$. As we saw above the formal neighborhood of $V$ can be identified with the completion $\hat{T}^*_V$ of the cotangent bundle on $V$. Let $x,y$ be local coordinates on $U:=\mathbb{C}\times\mathbb{C}\subset\mathbb{P}^1\times\mathbb{P}^1$. Then the derivations $\partial_x:=\partial/\partial x$ and $\partial_y:=\partial/\partial y$ are considered as (formal) functions on $\hat{T}^*_V$. One can check that $\mathcal{O}_{\hat{T}^*_V}(U)$ is isomorphic to $\mathbb{C}[x,y]\jump{\partial_x,\partial_y}$ and that the ring of global (formal) functions on $\hat{T}^*_V$ are generated by $\partial_x,x\partial_x,x^2\partial_x,\partial_y,y\partial_y$ and $y^2\partial_y$.

Note that $\hat{T}^*_V$ has a natural Poisson structure and Poisson brackets for local functions are given as follows:\\
$$\{\partial_x,x\}=\{\partial_y,y\}=1,\;\{\partial_x,\partial_y\}=\{\partial_x,y\}=\{\partial_y,x\}=\{x,y\}=0.$$
For other functions, the values of Poisson brackets can be computed using the linearity and the Leibniz rule.

We construct an explicit model of the formal neighborhood of $W$ intersecting with $\hat{T}^*_V$. Let $W_1,W_2$ be two copies of $\mathbb{C}\times\mathbb{P}^1$ and let $z_i$ and $w_i$ be local coordinates on $\mathbb{C}\times\mathbb{C}\subset\mathbb{C}\times\mathbb{P}^1$ for $i=1,2$. Then, similarly to above, the functions of $\hat{T}^*_{W_i}$ on this open subset is generated by $z_i,w_i,\partial_{z_i}:=\partial/\partial z_i$ and $\partial_{w_i}:=\partial/\partial_{w_i}$. We will glue $\hat{T}^*_{W_i}$'s to obtain a formal neighborhood of $W$. To do this, we first glue $\hat{T}^*_{W_i}$ with $\hat{T}^*_V$ along the curves $\{w_i=0\}$ and $s$.

Set $x'=1/x,\,y'=1/y,\,\partial_{x'}:=\partial/\partial x'$ and $\partial_{y'}:=\partial/\partial y'$. Then we have the relations $\partial_{x'}=-x^2\partial_x$ and $\partial_{y'}=-y^2\partial_y$. We identify the completion of $\hat{T}^*_{W_1}$ along $\{w_1=0\}$ and (an open formal subscheme of) the completion of $\hat{T}^*_V$ along $s$ via the following correspondence of functions:
$$z_1\leftrightarrow \frac{1}{2}(x+y),\;w_1\leftrightarrow -\frac{1}{2}(\partial_x-\partial_y),\;\partial_{z_1}\leftrightarrow \partial_x+\partial_y,\;\partial_{w_1}\leftrightarrow x-y.$$
Note that this correspondence preserves the natural Poisson structures. Similarly, we do this for $W_2$ via
$$z_2\leftrightarrow \frac{1}{2}(x'+y'),\;w_2\leftrightarrow-\frac{1}{2}(\partial_{x'}-\partial_{y'}),\;\partial_{z_2}\leftrightarrow \partial_{x'}+\partial_{y'},\;\partial_{w_2}\leftrightarrow x'-y'.$$

These also give identification of $\hat{T}^*_{W_1}$ and $\hat{T}^*_{W_2}$ along the open subschemes $\{z_i\ne0\}$ near the curves $\{w_i=0\}$. For example, $z_2$ is expressed using the coordinates of $W_1$ as follows:
$$z_2=\frac{1}{2}(x'+y')=\frac{x+y}{2xy}=\frac{4z_1}{(2z_1+\partial_{w_1})(2z_1-\partial_{w_1})}=\frac{1}{z_1}+\frac{1}{4z_1^3}\partial_{w_1}^2+\frac{1}{16z_1^5}\partial_{w_1}^4+\cdots.$$
One can check that this identification can be extended to the whole of the open subschemes of $\hat{T}^*_{W_i}$'s and that their underlying schemes $W_1$ and $W_2$ are glued to become $\Sigma_2$. Thus we obtain a formal Poisson scheme $\mathcal{W}$ whose underlying scheme is $W$. Since the isomorphism class of the completion $\hat{\mathcal{W}}_s$ is determined as the completion of $\hat{T}^*_V$ along $s$, the formal scheme $\mathcal{W}$ is isomorphic to the formal neighborhood of $W$ in the symplectic resolution of type $\mathbf{A}_3$ by Lemma \ref{lem:S2}. Note that the triviality assumption in Lemma \ref{lem:S2} can be directly checked though the second order deformation is nontrivial (see Remark \ref{rem3}).

Next we consider the ring $R$ of global (formal) functions on $\mathcal{W}$. Any $f\in R$ restricts to a global function of the formal neighborhood of $s$ in $\hat{T}^*_V$. It is not hard to check that any function on the formal neighborhood of $s$ extends to the whole of $\hat{T}^*_V$. Thus, $R$ is also the ring of global functions on the union $\mathcal{W}\cup \hat{T}^*_V$. By the formal function theorem \cite[Thm. 11.1]{Ha}, the singularity type for $\mathbf{A}_3$ is determined by $R$. Since the singularity has the unique symplectic resolution $Y$ (see Remark \ref{rem2}), the formal neighborhood $\hat{Y}_F$ is uniquely determined by $R$. Therefore, $\hat{Y}_F$ is isomorphic to the union of $\mathcal{W}$ and $\hat{T}^*_V$ which we have constructed.

\vspace{3mm}

\begin{rem}\label{rem3}
Unlike $V\cong\Sigma_0$, the formal neighborhood of $W\cong\Sigma_2$ is not isomorphic to the cotangent bundle. To show this, it suffices to show that the completion of $V$ along the diagonal $s$ is not isomorphic to the completion of the line bundle $\mathcal{O}_s(2)$ since the completion of $\hat{T}^*_W$ along $s$ is isomorphic to the completion of the normal bundle $\mathcal{O}_s(-2)^{\oplus2}\oplus\mathcal{O}_s(2)$.

Let $\mathcal{I}\subset\mathcal{O}_V$ be the ideal sheaf of $s$. Although one can construct algebra splittings of $\mathcal{O}_V/\mathcal{I}^n\to\mathcal{O}_V/\mathcal{I}$ for $n>0$, one sees that the sequence
$$0\to \mathcal{I}^2/\mathcal{I}^3\to \mathcal{I}/\mathcal{I}^3\to \mathcal{I}/\mathcal{I}^2\to0$$
does not split for any choice of the algebra splittings.

We can also prove the same thing by considering the ring of global functions on $\hat{T}^*_W$, which is not isomorphic to $R$ (see the next remark for the structure of $R$).
\qed
\end{rem}

\vspace{3mm}

\begin{rem}\label{rem4}
We can directly calculate the $\mathbb{C}$-algebra $R$ using the gluing of $\mathcal{W}$ and $\hat{T}^*_V$ as follows. Recall that the set of functions on $\hat{\mathcal{W}}_s$ is generated by
$f_1=\partial_x,f_2=x\partial_x,f_3=x^2\partial_x,f_4=\partial_y,f_5=y\partial_y$ and $f_6=y^2\partial_y$. In the coordinates of $W_1$, we can write
{\footnotesize
$$\begin{aligned}
&f_1=\frac{1}{2}\partial_{z_1}-w_1,\,f_2=(z_1+\frac{1}{2}\partial_{w_1})(\frac{1}{2}\partial_{z_1}-w_1),\,f_3=(z_1+\frac{1}{2}\partial_{w_1})^2(\frac{1}{2}\partial_{z_1}-w_1),\\
&f_4=\frac{1}{2}\partial_{z_1}+w_1,\,f_5=(z_1-\frac{1}{2}\partial_{w_1})(\frac{1}{2}\partial_{z_1}+w_1),\,f_6=(z_1-\frac{1}{2}\partial_{w_1})^2(\frac{1}{2}\partial_{z_1}+w_1).
\end{aligned}$$
}
To extend the functions on $\hat{\mathcal{W}}_s$ to the ones on $\mathcal{W}$, it suffices to extend to the ones on $\hat{T}^*_{W_1}$. Note that $f_i$'s themselves do not extend to $\hat{T}^*_{W_1}$ because of the terms $w_1,z_1w_1,z_1^2w_1$ and so on. We see that the following six functions are generators of $R$
{\footnotesize
$$\begin{aligned}
&a_1:=f_1+f_4=\partial_{z_1},\,a_2:=f_2+f_5=z_1\partial_{z_1}-w_1\partial_{w_1},\\
&a_3:=f_3+f_6=z_1^2\partial_{z_1}-2z_1w_1\partial_{w_1}+\frac{1}{4}\partial_{z_1}\partial_{w_1}^2,\,a_4:=f_1f_5-f_2f_6=(-\frac{1}{4}\partial_{z_1}^2+w_1^2)\partial_{w_1},\\
&a_5:=f_1f_6-f_3f_4=2z_1a_4,\,a_6:=f_2f_6-f_3f_5=(z_1^2-\frac{1}{4}\partial_{w_1}^2)a_4.
\end{aligned}$$
}
Note that $w_1\partial_{w_1}$ and $w_1^2\partial_{w_1}$ are global functions on $\hat{T}^*_{W_1}$. The relations of $a_i$'s should be the same as the ones given in section \ref{2} up to analytic coordinate change.
\qed
\end{rem}

\vspace{3mm}

We next consider the case where $\Gamma$ is of type $\mathbf{D}_n\,(n\ge4)$ and the fiber of the point $q_1$ (see Subsection \ref{2.2}). We first treat the case $\mathbf{D}_4$. The fiber $F$ in the symplectic resolution $Y$ is a union of the fiber $V\cup W$ for $\mathbf{A}_3$ and one more component $W'\cong\Sigma_2$. $W'$ intersects with $W$ along the negative section $s'$ of $W'$ which is disjoint from $s$. We construct an explicit model of formal neighborhood of $W'$ similarly to the case for $\mathbf{A}_3$. We use the same $W_i$'s and the coordinates as the case for $\mathbf{A}_3$. We set $v_i:=1/w_i$ and may assume that $s'\subset W$ is locally defined as $\{v_i=0\}$. Similarly to the case $\mathbf{A}_3$, we prepare the copies $W'_i\,(i=1,2)$ of $\mathbb{C}\times\mathbb{P}^1$ and let $z'_i$ and $w'_i$ be local coordinates on $\mathbb{C}\times\mathbb{C}\subset\mathbb{C}\times\mathbb{P}^1$ so that $s'\subset W'$ is locally defined as $\{w'_i=0\}$. We identify the completions of $W_i$ and $W'_i$ along $s'$ via the following correspondence:
$$z_i\leftrightarrow z'_i,\hspace{3mm}\partial_{z_i}\leftrightarrow\partial_{z'_i},\hspace{3mm}v_i\leftrightarrow\partial_{w'_i},\hspace{3mm}\partial_{v_i}\leftrightarrow -w'_i.$$
This gives identification of $\hat{T}^*_{W'_1}$ and $\hat{T}^*_{W'_2}$ along the open subschemes $\{z'_i\ne0\}$ near the curves $\{w'_i=0\}$. This extends globally and we obtain a formal Poisson scheme $\mathcal{W}'$ whose underlying scheme is $W'$. This is isomorphic to $\hat{Y}_{W'}$ by Lemma \ref{lem:S2}.

To use the same argument as for $\mathbf{A}_3$, we should show that any global function $f$ on the completion $\hat{\mathcal{W}}_{s'}$ extends to one on $\mathcal{W}$. To show this, it suffices to show that $f$ extends to the whole of $\hat{T}^*_{W_1}$. Assume that this were not true in order to deduce contradiction. $f$ can be written as a formal power series of $z_1,\partial_{z_1},v_1$ and $\partial_{v_1}$. Note that we can write $v_1=(\partial_x-\partial_y)^{-1}$ and $\partial_{v_1}=(\partial_x-\partial_y)^2(x-y)$. Since $f$ is assumed not to extend to $\hat{T}^*_{W'_1}$, the function $f$ should have $\partial_x-\partial_y$ as a pole when we write $f$ using the coordinates of $\hat{T}^*_V$. On the other hand, since $f$ is defined on the whole of $\hat{\mathcal{W}}_{s'}$, it is also written as a formal power series of $z_2,\partial_{z_2},v_2$ and $\partial_{v_2}$. Since these are written as
$$z_2=\frac{x+y}{2xy},\,\partial_{z_2}=x^2\partial_x+y^2\partial_y,\,v_2=-2(x^2\partial_x-y^2\partial_y)^{-1},\,\partial_{v_2}=\frac{(x-y)(x^2\partial_x-y^2\partial_y)^2}{-4xy},$$ 
the function $f$ cannot have poles along $\partial_x-\partial_y$, and hence contradiction. Therefore, we can show that $\hat{Y}_F$ is isomorphic, as a formal Poisson scheme, to the union of $\hat{T}^*_V\cup\mathcal{W}$ and $\mathcal{W}'$ which we have constructed. For higher $n\ge5$, by induction we can use the same argument as the case $n=4$.

\vspace{3mm}

\begin{rem}\label{rem5}
We can also directly calculate the $\mathbb{C}$-algebra $R_n$ of global functions on $\hat{Y}_F$ for $(\mathbf{D}_n,q_1)$ with $n\ge3$. Here we regard the ($n$=3)-case as the $\mathbf{A}_3$-case. Let $a_1,\dots,a_6$ be the functions on $\mathcal{W}\cup\hat{T}^*_V$ defined in Remark \ref{rem4}. Then $R_n$ is generated by the functions $a_1,a_2,a_3,b_{n,1},b_{n,2},b_{n,3}$ where $b_{n,i}$'s are defined by the following recurrence relations
$$b_{n+1,1}=\frac{1}{2}a_1 b_{n,2}-a_2 b_{n,1},\,b_{n+1,2}=a_1 b_{n,3}-a_3 b_{n,1},\,b_{n+1,3}=a_2 b_{n,3}-\frac{1}{2}a_3 b_{n,2}$$
with the initial terms $b_{3,1}=a_4,\,b_{3,2}=a_5$ and $b_{3,3}=a_6$. These generators have the following relation
$$a_1b_{n,3}-a_2b_{n,2}+a_3b_{n,1}=(a_1a_3-a_2^2)^{n-1}-b_{n,2}^2+4b_{n,1}b_{n,3}=0.$$

\vspace{-5mm}
\qed
\end{rem}

\vspace{3mm}

We have considered the gluing between $\Sigma_2$ and another component $\Sigma_0$ or $\Sigma_2$. Next we consider the gluing between two $\Sigma_0$'s. Let us consider the following situation. Let $Y$ be a 4-dimensional symplectic manifold and let $\pi:Y\to Y'$ be a contraction such that its exceptional locus $E$ is an irreducible divisor and that $\pi|_E$ is a $\mathbb{P}^1$-bundle whose image $S$ is a smooth surface. Assume that $S$ is the minimal resolution of an $\mathbf{A}_2$-singularity and let $C,C'\subset S$ be the irreducible exceptional curves. We also assume that both $V:=\pi^{-1}(C)$ and $W:=\pi^{-1}(C')$ are isomorphic to $\Sigma_0$. Moreover, we assume that there is a birational morphism from $Y$ which contracts the union $V\cup W$ to a point. In particular, there are divisors $E'$ and $E''$ of $Y$ such that\\
$\bullet$ $V=E\cap E'$ and $W=E\cap E''$, and\\
$\bullet$ there is a contraction $\psi:Y\to Z$ whose exceptional locus is $E'\cup E''$. 

We show that the isomorphism class of the formal neighborhood of $V\cup W$ is uniquely determined as the product (in the formal sense) of the minimal resolutions of $\mathbf{A}_1$ and $\mathbf{A}_2$ singularities. To show this, let $Y_1$ and $Y_2$ be two such symplectic manifolds and let $S_i$'s be the corresponding surfaces for $i=1,2$. Similarly, we define $V_i$, $W_i$, $E_i$ and so on. We fix an isomorphism of the unions $V_1\cup W_1\cong V_2\cup W_2$ that maps $V_1$ to $V_2$.

\begin{lem}\label{lem:sigma0}
Assume that we are in the situation as above. Also assume that we are given a Poisson isomorphism $\phi:(\hat{Y_1})_{V_1}\cong(\hat{Y_2})_{V_2}$ satisfying the following conditions:\\
$\bullet$ the restriction $\phi|_{V_1}$ is the given isomorphism $V_1\cong V_2$,\\
$\bullet$ $\phi$ induces an isomorphism $(\hat{S_1})_{C_1}\cong(\hat{S_2})_{C_2}$ which can extend to an isomorphism $(\hat{S_1})_{C_1\cup C'_1} \cong (\hat{S_2})_{C_2\cup C'_2}$ whose restriction to $C_1\cup C'_1$ is induced from $V_1\cup W_1\cong V_2\cup W_2$.

Then $\phi$ extends to an isomorphism $(\hat{Y_1})_{V_1\cup W_1}\cong(\hat{Y_2})_{V_2\cup W_2}$ whose restriction $V_1\cup W_1\cong V_2\cup W_2$ is the fixed one.
\end{lem}

{\em Proof.} Let $E'_i$ and $E''_i$ be the divisors of $Y_i$ such that $V_i=E_i\cap E'_i$ and $W_i=E_i\cap E''_i$, and let $\psi_i:Y_i\to Z_i$ be the contraction map of $E'_i\cup E''_i$ for $i=1,2$. Then $l_i:=\psi_i(V_i\cup W_i)\cong \mathbb{P}^1$ is a ($-2$)-curve in the surface $S'_i:=\psi_i(E'_i\cup E''_i)$. Thus the (local) functions on $(\hat{S'_i})_{l_i}$ are generated by $\bar{x}_i$ and $\partial_{\bar{x}_i}$ where $\bar{x}_i$ is the local coordinate of $l_i\cong\mathbb{P}^1$ defined on an open subset $U_i\cong \mathbb{C}$ of $l_i$. One can check that the Poisson isomorphism $\phi$ preserves $E'_i$'s and induces a Poisson isomorphism of $(\hat{S'_i})_{l_i}$'s. We take $\bar{x}_i$'s and $\partial_{\bar{x}_i}$'s so that $\bar{x}_1$ and $\partial_{\bar{x}_1}$ are mapped to $\bar{x}_2$ and $\partial_{\bar{x}_2}$ respectively by the isomorphism induced by $\phi$.

The function $\partial_{\bar{x}_1}$ can be considered as a global function on $E'_1\cup E''_1$ which vanishes on $V_1\cup W_1$, and we can take a global function $\partial_{x_1}$ on $(\hat{Y_1})_{V_1\cup W_1}$ which is a lift of $\partial_{\bar{x}_1}$. This is done by regarding $\partial_{\bar{x}_1}$ as a function on $\gamma((\hat{S'_1})_{l_1})$ and taking a lift on $\mathcal{Y}$ where $\gamma:(\hat{Y_1})_{V_1\cup W_1}\to \mathcal{Y}$ is the contraction of $V_1\cup W_1$. Note that $\mathcal{Y}$ is an affine formal scheme. We may assume that $\partial_{x_1}|_{E_1}=0$ by taking a lift on $\mathcal{Y}$ so that its zero locus contains $\gamma(\psi_1(E_1))$. This is possible since $\partial_{\bar{x}_1}$ is zero at the point $\gamma(l_1)=\gamma((\hat{S'_1})_{l_1})\cap \gamma(\psi_1(E_1))$. One sees that $\partial_{x_1}$ defines a Cartier divisor $E_1+2D_1$ for some $D_1$ on the opposite open subset $U'_1$ where $\bar{x}'_1:=1/\bar{x}_1$ is defined. Then we can take a lift $x'_1$ of $\bar{x}'_1$ as a defining function of $D_1$ on $U'_1$. Similarly we can take a lift $x_1$ of $\bar{x}_1$. We may assume that the lifts satisfy $\{\partial_{x_1},x_1\}=1$ by multiplying $\partial_{x_1}$ by a (unique) global invertible function.

Since $\phi$ induces an isomorphism $(\hat{S_1})_{C_1\cup C'_1}\cong(\hat{S_2})_{C_2\cup C'_2}$ by assumption, it preserves the exceptional divisors $E_1$ and $E_2$. This implies that $\partial_{x_2}:=\phi(\partial_{x_1})$ defines a global function on $(\hat{Y_2})_{V_2\cup W_2}$ which is a lift of $\partial_{\bar{x}_2}$. Similarly, the lifts of $\bar{x}_1$ and $\bar{x}'_1$ are mapped by $\phi$ to the lifts of $\bar{x}_2$ and $\bar{x}'_2:=1/\bar{x}_2$ respectively.

Let $R_{V_i\cup W_i}$ be the ring of global functions on $(\hat{Y_i})_{V_i\cup W_i}$ and let $R'_{V_i\cup W_i}$ be the subring of $R_{V_i\cup W_i}$ consisting of functions whose Poisson brackets with $x_i$ and $\partial_i$ are 0. Then one sees that, for each global function on $(\hat{S_i})_{C_i\cup C'_i}$, there is a unique lift which is in $R'_{V_i\cup W_i}$, and thus $R'_{V_i\cup W_i}$ is isomorphic to the ring of global functions on $(\hat{S_i})_{C_i\cup C'_i}$. This follows from the fact that, for any Poisson automorphism $\theta$ of $\hat{T}_{\Sigma_0}^*$ which is the identity on $\Sigma_0$, the values $\theta(y)$ and $\theta(\partial_y)$ are uniquely determined by $\theta(x)$ (and $\theta(\partial_x)$) where $x$ and $y$ are local coordinates of $\Sigma_0=\mathbb{P}^1\times\mathbb{P}^1$.

Let $R_{S'_i}$ be the ring of global functions on $(\hat{S'_i})_{l_i}$. We have a natural homomorphism $R'_{V_i\cup W_i}\hat{\otimes}R_{S'_i}\to R_{V_i\cup W_i}$ (cf. the proof of \cite[Prop. 3.3]{K}) where $\hat{\otimes}$ denotes the completed tensor product. This is an isomorphism by Lemma \ref{lem:glue} since we already know that its restrictions to $(\hat{Y_i})_{V_i}$ and $(\hat{Y_i})_{W_i}$ are isomorphisms. Therefore, $\phi$ gives an isomorphism $(\hat{Y_1})_{V_1\cup W_1}\cong(\hat{Y_2})_{V_2\cup W_2}$, which is an extension of $\phi$ by construction.
\qed

\vspace{3mm}

\begin{rem}\label{rem6}
From the proof of the lemma one sees that, for given Poisson isomorphisms $(\hat{S_1})_{C_1\cup C'_1}\cong(\hat{S_2})_{C_2\cup C'_2}$ and $(\hat{S'_1})_{l_1}\cong(\hat{S'_2})_{l_2}$, there is a Poisson isomorphism $(\hat{Y_1})_{V_1\cup W_1}\cong(\hat{Y_2})_{V_2\cup W_2}$ that induces the isomorphisms of the surfaces. Note that this isomorphism is not unique and depends on the choice of lifts of  functions. For $\hat{T}_{\Sigma_0}^*$, we have the following simple example. Consider the nontrivial Poisson automorphism $\theta$ of $\hat{T}_{\Sigma_0}^*$ defined by
$$\theta(x)=x+\partial_y,\,\theta(y)=y+\partial_x,\,\theta(\partial_x)=\partial_x,\,\theta(\partial_y)=\partial_y$$
with local coordinates $x$ and $y$ of $\Sigma_0$. Then $\theta$ induces the identity maps of the two copies of $\hat{T}_{\mathbb{P}^1}^*$. Note that the induced automorphism of the first (resp. the second) $\hat{T}_{\mathbb{P}^1}^*$ is obtained by substituting zero for $\partial_y$ (resp. $\partial_x$) in $\theta(x)$ and $\theta(\partial_x)$ (resp. $\theta(y)$ and $\theta(\partial_y)$).

One also sees that we can generalize this lemma to the case when $S_i$'s are the minimal resolutions of arbitrary ADE-singularities by repeating the same argument as above.
\qed
\end{rem}

\vspace{3mm}

This lemma enables us to treat the remaining cases: $\mathbf{A}_n\,(n\ge4),\,\mathbf{E}_n$ and $\mathbf{D}_n\,(n\ge5)$ for $q_2$ (see subsection \ref{2.2}). We only treat the $\mathbf{D}_5$ case since the other cases can be treated using the same methods.

Let $Y\to X$ be a symplectic resolution of a 4-dimensional symplectic variety $X$ and assume that $Y$ contains a 2-dimensional fiber $F'$ which is isomorphic to the fiber $F\subset \mathcal{H}$ of $q_2$ for $\mathbf{D}_5$ where $\mathcal{H}$ denotes the unique symplectic resolution of $\mathrm{Hilb}^2(\Gamma)$. We fix an isomorphism $F\cong F'$. We denote the irreducible component of $F'$ that corresponds to $P_{i,j}$ (resp. $Q_2$) by $P'_{i,j}$ (resp. $Q'_2$). Then the formal neighborhoods of $V:=P_{1,3}\cup P_{1,1}\cup Q_1$ in $\mathcal{H}$ (see Figure \ref{fig5}) and the corresponding union $V'$ in $Y$ are isomorphic  formal Poisson schemes as we already see.

We consider the exceptional divisor $E_2\subset\mathcal{H}$ defined in Subsection \ref{2.2} and let $E'_2\subset Y$ be the corresponding divisor. Then $E_2$ contains $P_{1,3}$ and $P_{1,4}$. The isomorphism of $\hat{\mathcal{H}}_V$ and $\hat{Y}_{V'}$ restricts to an isomorphism of $\hat{\mathcal{H}}_{P_{1,3}}$ and $\hat{Y}_{P'_{1,3}}$. It also induces an isomorphism $(\hat{S_1})_{C_1}\cong(\hat{S_2})_{C_2}$ (here we use the notation in Lemma \ref{lem:sigma0} with respect to $E_2,P_{1,3}$ and $P_{1,4}$). Note that $(\hat{S_1})_{C_1}$ can be identified with $\psi(E_3\cap \hat{\mathcal{H}}_{P_{1,1}})$ where $\psi$ is the contraction map of the divisor $E_3$ in $\hat{\mathcal{H}}_V$ where $E_3$ is defined in Subsection \ref{2.2}. We may assume that the isomorphism $(\hat{S_1})_{C_1}\cong(\hat{S_2})_{C_2}$ satisfies the second condition of Lemma \ref{lem:sigma0} since any Poisson automorphism of $(\hat{S_1})_{C_1}$ which is the identity on $C_1$ lifts to $\hat{\mathcal{H}}_V$ (see Remark \ref{rem7}). Therefore, we can extend the isomorphism $\hat{\mathcal{H}}_V\cong\hat{Y}_{V'}$ to the formal neighborhoods of $P_{1,4}$ and $P'_{1,4}$.

\begin{rem}\label{rem7}
The fact that the automorphisms of $(\hat{S_1})_{C_1}$ lift to $\hat{\mathcal{H}}_V$ is shown as follows. Let $t$ be a local coordinate of $C_1\cong\mathbb{P}^1$. Since $S_1$ is symplectic, $(\hat{S_1})_{C_1}$ is isomorphic to the completion of the cotangent bundle of $C_1$. Thus, the local coordinate ring of $(\hat{S_1})_{C_1}$ is generated by $t$ and $\partial_t:=\partial/\partial t$. Then, for any Poisson automorphism $\theta$ of $(\hat{S_1})_{C_1}$ which is the identity on $C_1$, we can write
$$\theta(t)=t+\sum_{a\ge0,b\ge1}c_{a,b}t^a\partial_t^b\hspace{5mm},\hspace{5mm}c_{a,b}\in\mathbb{C}.$$
The value $\theta(\partial_t)$ is uniquely determined from $\theta(t)$ using the equality $\{\theta(\partial_t),\theta(t)\}=1$. Let $\theta_{0,n},\theta_1$ and $\theta_2$ be the automorphisms defined by $\theta_{0,n}(t)=t+c_{0,n}\partial_t^n\,(n\ge1),\,\theta_1(t)=t+c_{1,1}t\partial_t$ and $\theta_2(t)=t+c_{2,1}t^2\partial_t$, and set $\theta_3:=\iota\circ\theta_1\circ\iota$ and $\theta_4:=\iota\circ\theta_{0,1}\circ\iota$ where $\iota$ is the involution defined by $\iota(t)=1/t$ and $\iota(\partial_t)=\partial_{1/t}$. One sees that any Poisson automorphism of $(\hat{S_1})_{C_1}$ is constructed successively from lower order terms with respect to $b$ by composing these automorphisms. (Note that taking the commutator $\theta\circ\theta'\circ\theta^{-1}\circ\theta'^{-1}$ for automorphisms $\theta,\theta'$ increases the minimal value of $b$ for which $c_{a,b}\ne0$.) For these generators, one can construct their liftings directly. For example, the lift $\bar{\theta}_{0,n}$ of $\theta_{0,n}$ is expressed with the coordinates of $P_{1,3}\cong\hat{T}^*_{\Sigma_0}$ as follows:
$$\bar{\theta}_{0,n}(x)=x+c_{0,n}(\partial_x+\partial_y)^n,\,\bar{\theta}_{0,n}(y)=y+c_{0,n}(\partial_x+\partial_y)^n,\,\bar{\theta}_{0,n}(\partial_x)=\partial_x,\,\bar{\theta}_{0,n}(\partial_y)=\partial_y.$$
One can check that this defines a Poisson automorphism of $\hat{\mathcal{H}}_V$ which is the identity on $V$ (see the construction of $\mathcal{W}$ and $\mathcal{W}'$). 
\qed
\end{rem}

\vspace{3mm}

Let $l_1,l_2\subset P_{3,3}$ be the intersections of $P_{3,3}$ with $P_{1,1}$ and $P_{1,4}$ respectively. Note that $l_1$ and $l_2$ are rulings of $P_{3,3}$ meeting at one point. We define $l'_1,l'_2\subset P'_{3,3}$ similarly. Then we have an isomorphism of $\hat{\mathcal{H}}_{l_1\cup l_2}$ and $\hat{Y}_{l'_1\cup l'_2}$. One can check that this extends to an isomorphism of $\hat{\mathcal{H}}_{P_{3,3}}$ and $\hat{Y}_{P'_{3,3}}$. Similarly, we obtain an isomorphism of $\hat{\mathcal{H}}_{P_{1,2}}$ and $\hat{Y}_{P'_{1,2}}$ from the two isomorphisms for $P_{1,3}$ and $P_{3,3}$.

By repeating such processes, we can extend the isomorphism to the whole $F$ and $F'$. The same methods can apply to the other cases and we obtain the following theorem.

\begin{thm}\label{thm:main}
Let $X$ be a 4-dimensional symplectic variety and $\pi:Y\to X$ a symplectic resolution. Assume the fiber $\pi^{-1}(x)$ of $x\in X$ is isomorphic to a 2-dimensional fiber $\psi^{-1}(q)$ of the unique symplectic resolution $\psi:\mathcal{H}\to\mathrm{Hilb}^2(\Gamma)$ where $\Gamma$ is a surface with one ADE-singularity (cf. Section \ref{2}). Then the analytic germs $(X,x)$ and $(\mathrm{Hilb}^2(\Gamma),q)$ are isomorphic.
\end{thm}

{\em Proof.} We have shown that the 2-dimensional fibers have isomorphic formal neighborhoods. This implies that the singularity type of $\mathrm{Hilb}^2(\Gamma)$ at its 0-dimensional symplectic leaf is uniquely determined by the formal function theorem \cite[Thm. 11.1]{Ha}.
\qed

\section{Slodowy slices of nilpotent cones to sub-subregular orbits}\label{4}

In this section we introduce some kinds of Slodowy slices in simple Lie algebras $\mathfrak{g}$ of simply-laced types and study their symplectic resolutions.

Let $\mathfrak{g}$ be a finite-dimensional complex simple Lie algebra. It is well-known that such $\mathfrak{g}$ is completely classified (see e.g. \cite{Hu}). $\mathfrak{g}$ defines its adjoint group $G$, which linearly acts on $\mathfrak{g}$. The orbit of a nilpotent element in $\mathfrak{g}$ by $G$ is called a {\em nilpotent (adjoint) orbit}. Every nilpotent orbit $O\subset\mathfrak{g}$ is a nonsingular algebraic variety and admits a holomorphic symplectic form $\omega_{KK}$ called the {\em Kostant-Killirov form}. See \cite{CM} for generalities about nilpotent orbits.

\begin{prop}\cite{P}
The normalization of the closure $\overline{O}$ of a nilpotent orbit $O\subset\mathfrak{g}$ is a symplectic variety.
\end{prop}

Among the nilpotent orbits in $\mathfrak{g}$, which are in fact finite, there is the biggest one $O_{\mathrm{reg}}$ in the sense that its closure $\mathcal{N}$ (called the {\em nilpotent cone}) contains all the other nilpotent orbits. This orbit is called the {\em regular orbit}. Similarly, there is the second biggest orbit $O_{\mathrm{subreg}}$ called the {\em subregular orbit}. However, $\overline{O}_{\mathrm{subreg}}\setminus O_{\mathrm{subreg}}$ can be reducible in general. Therefore, there are non-unique third biggest orbits and we call them {\em sub-subregular orbits}.

Next we introduce Slodowy slices. For an element $x$ of a nilpotent orbit $O$, by Jacobson-Morozov theorem, there are two elements $y,h\in\mathfrak{g}$ such that
$$[x,y]=h,\,[h,x]=2x,\text{ and }[h,y]=-2y$$
(such $x,y$ and $h$ are called an $\mathfrak{sl}_2$-triple). Then we define an affine subspace of $\mathfrak{g}$ as
$$\mathcal{S}_x=x+\mathrm{Ker}\,(\mathrm{ad}\,y)$$
where $\mathrm{ad}\,y:\mathfrak{g}\to\mathfrak{g}$ is a $\mathbb{C}$-linear map defined by $v\mapsto [y,v]$. In fact the isomorphism class of the intersection $\mathcal{S}_x\cap\mathcal{N}$ does not depend on the choice of $x,y$ and $h$. Thus we denote it simply by $\mathcal{S}$ and call it the {\em Slodowy slice} (of the nilpotent cone to the nilpotent orbit $O$). The Slodowy slice is a symplectic variety with a symplectic form inherited from $\mathcal{N}$. From now on we consider the case when $O$ is a sub-subregular orbit.

\vspace{3mm}

{\bf Example}  The case $\mathfrak{g}=\mathfrak{sl}_4$.\\
We can take an $\mathfrak{sl}_2$-triple as
{\footnotesize
$$x=\begin{pmatrix}
0&1&0&0\\
0&0&0&0\\
0&0&0&1\\
0&0&0&0
\end{pmatrix},\;
y=\begin{pmatrix}
0&0&0&0\\
1&0&0&0\\
0&0&0&0\\
0&0&1&0
\end{pmatrix},\;
h=\begin{pmatrix}
1&0&0&0\\
0&-1&0&0\\
0&0&1&0\\
0&0&0&-1
\end{pmatrix}.$$}
Then
{\footnotesize
$$\mathcal{S}_x=
\left\{ \begin{pmatrix}
t_1&1&s_1&0\\
t_2&t_1&s_2&s_1\\
s_3&0&-t_1&1\\
s_4&s_3&u&-t_1
\end{pmatrix}\middle|\; t_1,t_2,s_1,s_2,s_3,s_4,u\in\mathbb{C}\right\}.$$}

By calculating the characteristic polynomial, we know that the Slodowy slice $\mathcal{S}_x\cap\mathcal{N}$ is defined by the following 2 equations in $\mathbb{C}[t_1, s_1, s_2, s_3, s_4, u]$ ($t_2$ is eliminated):\\
{\footnotesize
$$2t_1^3 +2t_1s_1s_3-s_2s_3-s_1s_4+2t_1u,$$
$$6t_1^2s_1s_3+6s_1^2s_3^2+3t_1s_2s_3+3t_1s_1s_4-10t_1^2u-4s_1s_3u- 2s_2s_4-2u^2.$$}
As we will see below, the analytic germ of the Slodowy slice at the origin is isomorphic to that of $U_5\subset\mathrm{Hilb}^2(\Gamma_3)$ at the origin (see Subsection \ref{2.1}).
\qed

\vspace{3mm}

It is well known that the nilpotent cone $\mathcal{N}$ has a canonical resolution $\pi:\tilde{\mathcal{N}}\to\mathcal{N}$ called the Springer resolution. Here $\tilde{\mathcal{N}}$ is the cotangent bundle of a certain rational homogeneous space and thus it has a natural symplectic structure. $\pi$ is a symplectic resolution with respect to this symplectic structure.

As for the Slodowy slice $\mathcal{S}$, its symplectic resolution is obtained as just the restriction of the Springer resolution to $\mathcal{S}$. From now on we concentrate on the cases where $\mathfrak{g}$ is of type $\mathbf{A}_n\,(n\ge3),\mathbf{D}_n\,(n\ge4)$ or $\mathbf{E}_n\,(n=6,7,8)$ and study the geometry of the symplectic resolution $p:\tilde{\mathcal{S}}\to\mathcal{S}$.

For $\mathbf{A}_n$ and $\mathbf{E}_n$, the choice of the sub-subregular orbit is unique. For $\mathbf{D}_4$, we have three choices which are symmetric (thus give the same singularity types). For $\mathbf{D}_n\,(n\ge5)$, we can choose two sub-subregular orbits: the ADE-type of one in the 2-dimensional symplectic leaf $\mathcal{S}_x\cap\overline{O}_{\mathrm{subreg}}$ is $\mathbf{A}_1$ and the other $\mathbf{D}_{n-2}$. The main theorem in this section is the following.

\begin{thm}\label{thm:slice}
Let $\Gamma$ be the singular surface as in Section \ref{2} whose type is the same as $\mathfrak{g}$. Then there is an isomorphism of germs
$$(\mathrm{Hilb}^2(\Gamma),q)\cong(\mathcal{S},x)$$
where $q$ is a 0-dimensional symplectic leaf of $\mathrm{Hilb}^2(\Gamma)$, which is unique when $\mathfrak{g}$ is of type $\mathbf{A}_n$ or $\mathbf{E}_n$. When $\mathfrak{g}$ is of type $\mathbf{D}_n$, $q$ is chosen so that the types of $q$ and $x$ are the same in the 2-dimensional symplectic leaves.
\end{thm}

To prove this theorem, we focus on fibers of the symplectic resolution as in the previous sections. The famous result of Brieskorn-Slodowy theory states that the fiber of an element of the subregular orbit is the Dynkin tree of $\mathbb{P}^1$'s of the same type as $\mathfrak{g}$.

The central fiber $p^{-1}(x)$ is more complicated. In \cite{Lo} Lorist explicitly described the irreducible components of $p^{-1}(x)$. By comparing his result with ours obtained in Section \ref{2}, we see that the central fibers of the symplectic resolutions of $\mathrm{Hilb}^2(\Gamma)$ and $\mathcal{S}$ are isomorphic. The description by Lorist also shows that the fibers for $\mathbf{E}_6$ and $\mathbf{E}_7$ have embeddings into the fiber for $\mathbf{E}_8$.

\vspace{3mm}

{\em Proof of Theorem \ref{thm:slice}.}\\
Since the central fibers of the symplectic resolutions of $\mathrm{Hilb}^2(\Gamma)$ and $\mathcal{S}$ are isomorphic, the claim follows from Theorem \ref{thm:main}.

\qed

\vspace{0.2cm}

\begin{center}
Department of Mathematics, Graduate School of Science, Kyoto University, Japan 

ryo-yama@math.kyoto-u.ac.jp
\end{center}

\end{document}